 \documentstyle [12pt,twoside, epsf]{article}
 \newcommand{\ZZ}{\mbox{$Z\!\!\! Z\!$}} 
 \newcommand{\QQ}{\mbox{$Q\!\!\!\! I$}} 

 \let\\=\cr
 
 \def\Chi{\hbox{\raise0.5ex\hbox{$\chi$}}}

 \newtheorem{th}{Theorem}
 \newtheorem{lem}{Lemma}
 
 \newtheorem{prop}{Proposition}

 \newtheorem{defn}{Definition}

 \newtheorem{rem}{Remark}

 \def\picill#1by#2(#3){\epsffile{#3}}

\begin{document}
 \pagestyle{myheadings}

  \markboth{{\sc Kauffman \& Lambropoulou}}{{\sc  From Tangle Fractions to DNA}}

 \title{From Tangle Fractions to DNA}

 \author{Louis H. Kauffman and Sofia Lambropoulou}

\date{}

 \maketitle


\begin{abstract}
This paper draws a line from the elements of tangle fractions to the tangle model of DNA recombination.
In the process, we sketch the classification of rational tangles, unoriented and oriented rational knots and the 
application of these subjects to DNA recombination.
\end{abstract}

\section{Introduction}

Rational knots and links
are a class of alternating links of one or two
unknotted components, and they are
the easiest knots to make (also for Nature!). 
The first twenty five knots,
except for $8_5$,
are rational.  Furthermore all knots and links up to ten crossings are either rational
or are obtained
from rational knots by insertion operations on certain simple graphs.
Rational knots are also known in the literature as four-plats,
Viergeflechte and $2$-bridge knots.  
The lens spaces arise as $2$-fold branched coverings along rational knots.

A rational tangle is the result
of consecutive
twists on  neighbouring endpoints of two trivial arcs, see Definition 1.
Rational knots are
obtained by taking numerator closures of rational tangles (see Figure 19),
which form a basis
for their classification. Rational knots and rational tangles are of
fundamental importance in
the study of DNA recombination.  
Rational knots and links were first considered in \cite{Rd2} and \cite{BS}.
Treatments  of various  aspects of rational
knots and rational
tangles can be found in \cite{Bl},\cite{C1}, \cite{Sie}, \cite{BZ}, \cite{R}, \cite{GK2}, \cite{Kaw},
\cite{Li}, \cite{M}.  
A rational tangle  is associated in a canonical manner with a unique, reduced
rational number or $\infty,$ called {\it the fraction} of the tangle.
Rational tangles are
 classified by their fractions by means of the following theorem:

 \begin{th}[Conway, 1970]{ \ Two rational tangles are isotopic if and only
if
 they have the same fraction.
  } \end{th}

John H. Conway  \cite{C1} introduced the notion of tangle and
defined the fraction of a
rational tangle using the continued fraction form of the tangle and the 
Alexander polynomial of knots.  Via the Alexander polynomial, the fraction is 
defined for the larger class of all $2$-tangles. In this paper we are interested in different definitions of the 
fraction, and we give a self-contained exposition of the construction of the invariant fraction for arbitrary 
$2$-tangles from the bracket polynomial \cite{K1}. The tangle fraction is a key ingredient in both the classification of
rational knots and in the applications of knot theory to DNA.  
Proofs of Theorem~1 can be found in
\cite{Mo}, \cite{BZ} p.196, \cite{GK2} and \cite{KL1}. 
\smallbreak

More than one rational tangle can yield the same or isotopic
rational knots and the equivalence relation between the rational tangles is
reflected in an
arithmetic equivalence of their corresponding fractions. This is marked by a
theorem
due originally to Schubert \cite{Sch} and reformulated by  Conway \cite{C1}
in  terms of rational tangles.

 \begin{th}[Schubert, 1956]{ \ Suppose that rational tangles with fractions
$\frac{p}{q}$ and
 $\frac{p'}{q'}$ are given ($p$ and $q$ are relatively prime.
 Similarly for $p'$ and $q'$.) If  $K(\frac{p}{q})$ and $K(\frac{p'}{q'})$
denote the
 corresponding rational knots obtained by taking numerator closures of
these tangles, then
$K(\frac{p}{q})$ and $K(\frac{p'}{q'})$ are topologically equivalent if and
only if
 \begin{enumerate}
 \item $p=p'$ and
 \item either $q\equiv q'(mod \, p)$  \ or \ $qq'\equiv 1(mod \, p).$
 \end{enumerate}
  } \end{th}

This classic theorem  \cite{Sch} was originally proved by
using an observation of Seifert that
the $2$-fold
branched covering spaces of $S^3$ along $K(\frac{p}{q})$
and $K(\frac{p'}{q'})$ are lens spaces,
and invoking the results of Reidemeister \cite{Rd3} on the classification of lens
spaces. Another proof
using covering spaces has been given by  Burde in \cite{Bu}. Schubert also
extended this theorem to
the case of oriented rational knots and links described as
$2$-bridge links:

 \begin{th}[Schubert, 1956]{ \ Suppose that orientation-compatible rational
tangles with
 fractions $\frac{p}{q}$  and $\frac{p'}{q'}$ are given with $q$ and $q'$
odd. ($p$ and $q$ are
 relatively prime. Similarly for $p'$ and $q'$.) If
$K(\frac{p}{q})$ and
$K(\frac{p'}{q'})$ denote the corresponding rational knots obtained by
taking numerator
closures of these tangles, then $K(\frac{p}{q})$ and $K(\frac{p'}{q'})$ are
topologically
equivalent if and only if

 \begin{enumerate}
 \item $p=p'$ and
 \item either $q\equiv q'(mod \, 2p)$  \ or \ $qq'\equiv 1(mod \, 2p).$

 \end{enumerate}
  } \end{th}

In \cite{KL2} we give the first combinatorial proofs of Theorem~2 and Theorem~3. 
In this paper we sketch the proofs
in \cite{KL1} and \cite{KL2} of the above three
theorems and we give the key examples that are behind all of our proofs.
We also give
some  applications of
Theorems~2 and 3 using our methods.
\bigbreak

 The paper is organized as follows. In Section~2 we introduce $2$-tangles and rational tangles,
Reideimeister moves,
 isotopies and operations. We give the definition of flyping, and state the (now proved) Tait flyping conjecture.
The Tait conjecture is used implicitly in our classification work.  In Section~3 we
introduce the continued fraction expression for rational tangles and its properties. We use the continued fraction 
expression for rational tangles to define their fractions. Then rational tangle diagrams are shown to be isotopic to
alternating diagrams. The alternating form is used to obtain a canonical form for rational tangles, and we obtain a proof
of Theorem~1. 

Section~4 discusses alternate definitions of the tangle fraction. We begin with a self-contained exposition of
the bracket polynomial for knots, links and tangles. Using the bracket polynomial we define a fraction $F(T)$ 
for arbitrary $2$-tangles and show that it has a list of properties that are sufficient to prove that for $T$ rational, $F(T)$ is
identical to  the continued fraction value of $T$, as defined in Section~3. The next part of Section~4 gives a
different definition of the fraction of a rational tangle, based on coloring the tangle arcs with integers. This definition is
restricted to rational tangles and those tangles that are obtained from them by tangle-arithmetic operations, but it is truly
elementary,  depending just on a little algebra and the properties of the Reidemeister moves. Finally, we sketch yet another
definition of the  fraction for $2$-tangles that shows it to be the value of the conductance of an electrical network associated
with the tangle. 

Section~5 contains a description of our approach to the proof of Theorem~2, the classification of unoriented rational knots and
links. The key to this approach is enumerating the different rational tangles whose numerator closure is a given unoriented
rational knot or link, and confirming that the corresponding fractions of these tangles satisfy the arithmetic relations of the
Theorem.  Section~6 sketches the classification of  rational knots and links that are isotopic to their mirror images. Such links
are all closures of  palindromic continued fraction forms of even length. Section~7 describes our proof of Theorem~3, the
classification of oriented rational knots. The statement of Theorem~3 differs from the statement of Theorem~2 in the use of
integers modulo $2p$ rather than $p.$ We see how this difference arises in relation to matching orientations on
tangles.  This section also includes an
explanation of the fact that fractions with even numerators correspond to rational links of two components, while fractions with
odd numerators correspond to single component rational knots (the denominators are odd in both cases). Section~8 discusses strongly
invertible rational knots and links. These correspond to palindromic continued fractions of odd length. 

Section~9 is an
introduction to the tangle model for DNA recombination. The classification of the rational knots
and links, and the use of the tangle fractions is the basic topology behind the tangle model for DNA recombination.
 We indicate how problems in this model are reduced to
properties of rational knots, links and tangles, and we show how a finite number of observations of successive DNA recombination
can pinpoint the recombination mechanism.

 \section{$2$-Tangles and Rational Tangles}

 Throughout this paper we will be working with {\it $2$-tangles}. The theory of
 tangles was
 discovered  by John Conway \cite{C1} in his work on enumerating and
classifying knots.
 A $2$-tangle is an embedding of two arcs
(homeomorphic to
 the interval [0,1]) and circles into a three-dimensional ball $B^3$  standardly embedded in Euclidean three-space $S^3$,
 such that the endpoints of the arcs go to a specific set of four  points
on the surface of the
 ball, so that the circles
 and the interiors of the arcs are embedded in the interior of the ball.  
The left-hand side of Figure 1 illustrates
a $2$-tangle. Finally,
 a $2$-tangle is  {\it oriented} if we assign orientations to each arc
and each circle. Without loss of
 generality, the four endpoints of a $2$-tangle can be arranged on a great
circle on the boundary of the ball. One can then  define a {\it diagram} of a
 $2$-tangle to be a regular projection of the tangle on the plane of
 this great circle.  In illustrations we may replace this circle by a box.  
\bigbreak

$$ \picill3.5inby4.5in(D1) $$

 \begin{center}
 {\bf Figure 1 - A $2$-tangle and a rational tangle }
 \end{center}
\vspace{3mm}

   The simplest possible $2$-tangles comprise
two unlinked arcs either
horizontal or vertical. These are the {\it trivial
 tangles}, denoted  $[0]$ and $[\infty]$ tangles respectively, see Figure 2.
\bigbreak

$$ \picill3.5inby4.5in(D2) $$

\begin{center}
{\bf Figure 2 - The trivial tangles $[0]$ and $[\infty]$ }
\end{center}
\vspace{3mm}

 \begin{defn}{\rm \ A
 $2$-tangle is {\it rational}  if it can be obtained
by applying a finite number of consecutive twists of neighbouring endpoints to
the elementary tangles
$[0]$ or $[\infty].$} 
\end{defn}

The simplest rational tangles are the $[0]$, the $[\infty]$, the $[+1]$ and the
$[-1]$ tangles, as illustrated in Figure 3, while the next
simplest ones are:

 \begin{itemize}
 \item[(i)] \ The  {\it integer tangles}, denoted by $[n],$ made of $n$
 horizontal twists, $n \in
 \ZZ.$
 
 \item[(ii)] \ The {\it vertical tangles},  denoted by $\frac{1}{[n]},$
 made of $n$ vertical twists, $n~\in~\ZZ.$  These are the inverses of the
integer tangles, see
Figure 3.
 This terminology will be clear soon.
 \end{itemize}
\bigbreak

 Examples of rational tangles are
illustrated in the right-hand
side of Figure 1 as well as in  Figures 8 and 17 below.
\bigbreak

$$ \picill4.5inby1.65in(D3) $$
 \begin{center}
 {\bf Figure 3 - The elementary rational tangles}
 \end{center}
 \vspace{3mm}

We study tangles up to 
{\it isotopy}. Two
$2$-tangles, $T, S$, in $B^3$ are said to be  {\it isotopic}, denoted by
$T \sim S$, if they have
identical configurations of their four endpoints in the boundary $S^2$ of the three-ball,
and there is an ambient
isotopy of $(B^3, T)$ to $(B^3, S)$ that is the identity on the boundary 
$(S^2,\partial T) =
(S^2,\partial S)$. An ambient isotopy can be imagined as a continuous
deformation of $B^3$ fixing the
four endpoints on the boundary sphere, and bringing one tangle to the other
without causing any
self-intersections.   
 \smallbreak

In terms of diagrams, Reidemeister \cite{Rd} proved that the
local moves on diagrams
illustrated in Figure 4 capture combinatorially the notion of ambient
isotopy of knots, links and
tangles in three-dimensional space. That is, if two diagrams represent
knots, links or tangles that
are isotopic, then the one diagram can be obtained from the other by a
sequence of {\it Reidemeister
moves}. In the case of tangles {\it the endpoints of the tangle remain
fixed} and all the moves occur
inside the tangle box.
\smallbreak

Two oriented $2$-tangles are are said to be  {\it oriented
isotopic} if there is an
isotopy between them that preserves the orientations of the corresponding
arcs and the corresponding
circles. The diagrams of two oriented isotopic tangles differ by a sequence
of {\it oriented Reidemeister moves}, i.e. Reidemeister moves with orientations
on the little arcs that remain
consistent during the moves. 
\bigbreak

$$ \picill3.85inby2.5in(D4)  $$

\begin{center}
{\bf Figure 4 - The Reidemeister moves}
\end{center}
\vspace{3mm}

From now on we will be thinking in terms of
tangle diagrams. Also, we will be referring to both knots and links whenever we say
`knots'.
\smallbreak

A {\it flype} is an isotopy move  applied on a $2$-subtangle of a larger tangle or knot as shown in Figure 5. A
 flype preserves the
 alternating structure of a diagram. Even more, flypes are the only isotopy
 moves needed in the
 statement of the celebrated Tait Conjecture for alternating knots, stating
 that {\it two alternating
 knots
  are isotopic if and only if any two corresponding diagrams on $S^2$ are
 related by a finite
 sequence of flypes.} This was posed by P.G.~Tait, \cite{Ta} in 1898 and
 was
 proved by
 W.~Menasco and M.~Thistlethwaite, \cite{MT} in 1993.
\bigbreak

$$ \picill2.95inby2in(D5) $$

\begin{center}
 {\bf Figure 5 - The flype moves}
 \end{center}
 \vspace{3mm}

 The class of $2$-tangles is closed under the operations of {\it addition}
 ($+$) and {\it multiplication}  ($*$) as illustrated  in  Figure 6. Adddition is
 accomplished
 by placing the tangles
 side-by-side and attaching the
 $NE$ strand of the left tangle to the $NW$ strand of  the right tangle,
 while attaching the $SE$
 strand of the left tangle to the $SW$ strand of the right tangle. The 
 product is accomplished by
 placing one tangle underneath the other and attaching the upper strands of
 the lower tangle to the
 lower strands of the upper tangle.
 \smallbreak

 The {\it mirror image} of  a tangle $T$ is
 denoted by $-T$ and it is obtained by switching all the crossings in $T.$
Another operation is {\it rotation} accomplished by turning the tangle counter-clockwise by
$90^{\circ}$ in the plane. The rotation of $T$ is denoted by $T^r.$
The {\it inverse} of a tangle $T$, denoted by
 $1/T,$ is defined to be  $-T^r.$ See Figure 6.
In general, the inversion or rotation  of a $2$-tangle is an order 4 operation.   
 Remarkably, for  rational tangles the
 inversion (rotation) is an order 2 operation.  
 It is for this reason that we denote the inverse of a $2$-tangle $T$ by
 $1/T$ or $T^{-1},$ and hence the rotate  of the tangle $T$ can be denoted by
 $-1/T = -T^{-1}.$  
\bigbreak

$$\vbox{\picill4.25inby1.82in(D6)  }$$
\begin{center}
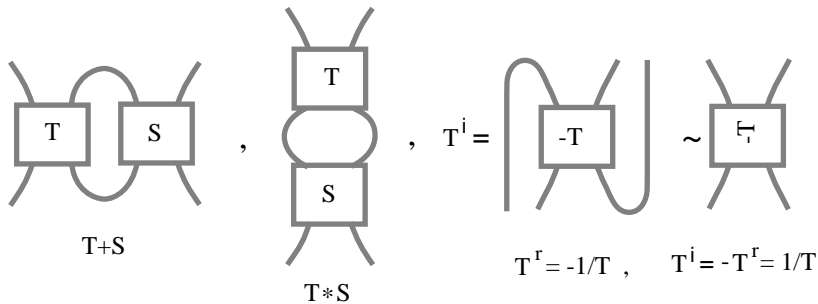
 {\bf Figure 6 - Addition, product and inversion of $2$-tangles }
 \end{center}
\vspace{3mm}

 We describe now another operation applied  on $2$-tangles, which
 turns
 out to be an isotopy on rational tangles. We
 say that
 $R^{hflip}$ is the {\it horizontal flip} of the tangle $R$ if $R^{hflip}$
 is
 obtained from
 $R$ by a $180^{\circ}$ rotation around a horizontal axis on the plane of $R$.
 Moreover, $R^{vflip}$
 is the  {\it vertical flip} of the $2$-tangle $R$ if
 $R^{vflip}$ is obtained from $R$ by a $180^{\circ}$ rotation around a vertical
 axis on the plane of $R$.
 See Figure 7 for illustrations.  Note that a flip switches the
 endpoints of the tangle and, in general, a flipped tangle is not isotopic
 to
 the original one.
 {\it It is a property of rational tangles that $T \sim
 T^{hflip}$ and $T \sim
 T^{vflip}$ for any rational tangle $T.$} This
 is obvious for the
 tangles $[n]$ and
 $\frac{1}{[n]}.$ The general proof crucially uses flypes, see \cite{KL1}.  
\bigbreak

 $$ \picill3inby1.8in(D7) $$
 \begin{center}
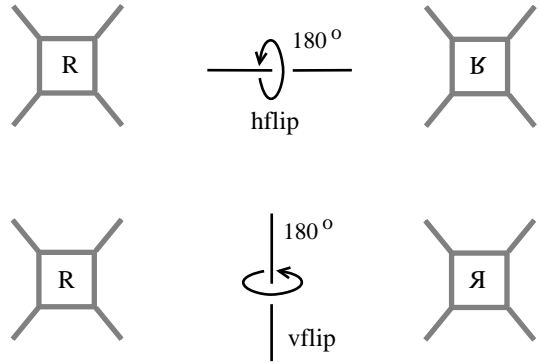
 {\bf Figure 7 - The horizontal and the vertical flip}
\end{center}
 \vspace{3mm}

 The above isotopies composed consecutively yield $T \sim
 (T^{-1})^{-1}
 = (T^{r})^{r}$ for any rational tangle $T.$  This says that inversion
 (rotation) is an operation of
 order 2 for rational tangles, so we
 can rotate the mirror image
 of $T$  by $90^{\circ}$ either counterclockwise  or clockwise to obtain $T^{-1}$.
 \bigbreak

Note
that the twists generating
the rational tangles could take place between the right, left, top or
bottom  endpoints of a
previously created rational tangle. Using flypes and flips inductively on subtangles
one can always
bring the twists to the right or bottom of the rational tangle.  We shall then say
that the rational tangle is
in {\it standard form}. Thus a rational tangle in standard form  is created
by consecutive
additions of the tangles $[\pm 1]$ {\it only on the right}  and multiplications
by the tangles $[\pm 1]$
{\it only at the bottom,} starting from the tangles $[0]$ or $[\infty].$ 
For example, Figure 1
 illustrates the tangle  $(([3]*\frac{1}{[-2]})+[2]),$ while  Figure 17
 illustrates the tangle  $(([3]*\frac{1}{[2]})+[2])$ in standard form.
Figure 8 illustrates addition on the right and multiplication on the bottom by elementary tangles.

\bigbreak

 $$ \picill3.1inby2.3in(D8) $$
\begin{center}
{\bf Figure 8 - Creating new rational tangles}
 \end{center}
\vspace{3mm}

 We also have the following
 {\it closing} operations, which yield two different knots: the {\it
 Numerator} of a $2$-tangle $T$,
 denoted by $N(T)$, obtained by joining with simple arcs the two
 upper endpoints and the two
 lower endpoints of $T,$ and the  {\it Denominator} of a $2$-tangle $T$,
 obtained
 by joining with simple arcs each pair of the corresponding top and bottom
 endpoints of $T$,
 denoted by $D(T)$. We have $N(T) = D(T^r)$ and $D(T) =
 N(T^r).$   We note that every
 knot or link can be regarded as the numerator closure of a $2$-tangle.
\bigbreak

 $$\vbox{\picill4.3inby1.37in(D9)  }$$
 \begin{center}
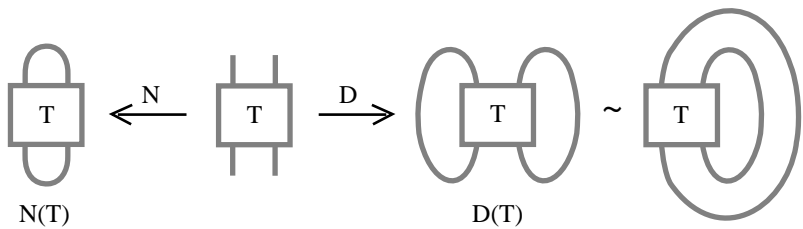
 {\bf Figure 9 - The numerator and denominator of a $2$-tangle}
 \end{center}
 \vspace{3mm}

 We obtain $D(T)$ from $N(T)$ by a
  $[0]-[\infty]$ interchange, as shown in Figure~10. This `transmutation' of
 the numerator to the
 denominator is a precursor to the tangle model of a recombination event in
 DNA, see Section 9.
 The $[0]-[\infty]$ interchange can be described algebraically by the
 equations:

 $$N(T) = N(T+[0]) \longrightarrow N(T+[\infty]) = D(T).$$
\bigbreak

$$\vbox{\picill4.25inby1.5in(D10)  }$$
 \begin{center}
 {\bf Figure 10 - The $[0]-[\infty]$ interchange } \end{center}
 \vspace{3mm}

 We will concentrate on the
 class of {\it rational knots and links}  arising from closing the
 rational tangles.
 Even though the sum/product of rational tangles is in
 general not rational,
  the numerator (denominator) closure of the sum/product of two rational
 tangles is still a
 rational knot. It  may happen that
 two rational tangles are not isotopic but have isotopic  numerators.  This
 is the basic idea
 behind the classification of rational knots, see Section 5.
 \bigbreak

\section{Continued Fractions and the Classification of Rational Tangles}

 In this section we assign to a rational tangle a fraction, and we explore
 the analogy between rational
 tangles and continued fractions. This analogy culminates in a common
 canonical form, which is used to
 deduce the classification of rational tangles.
\bigbreak

We first observe that multiplication of a rational tangle $T$ by
$\frac{1}{[n]}$
 may be obtained as addition of $[n]$ to the inverse $\frac{1}{T}$ followed
by inversion.
Indeed, we have:

 \begin{lem}{ \ The following tangle equation holds for any rational tangle
 $T.$

 $$ T* \frac{1}{[n]} = \frac{1}{[n] + \frac{1}{T}}. $$ Thus any rational tangle can be built by a series
of the following operations:  {\it Addition
 of $[\pm 1]$} and
 {\it Inversion}. 

  } \end{lem}

 \noindent {\em Proof.}  Observe that a $90^{\circ}$ clockwise rotation of
 $T* \frac{1}{[n]}$ produces $-[n] - \frac{1}{T}.$ Hence, from the above
 ${(T* \frac{1}{[n]})}^r = -[n]
 - \frac{1}{T},$ and thus $({T* \frac{1}{[n]}})^{-1} = [n] + \frac{1}{T}.$ So,
 taking inversions on both
 sides yields the tangle equation of the statement. $\hfill \Box $

 \begin{defn}{\rm \ A {\it continued fraction in integer tangles} is an
 algebraic description of a
 rational tangle via a continued fraction  built from the  tangles $[a_1],$
 $[a_2],$
 $\ldots, [a_n]$  with all numerators equal to~$1$, namely an expression of
 the type:

 \[ [[a_1],[a_2],\ldots ,[a_n]] \, := \, [a_1]+ \frac{1}{[a_2]+\cdots +
 \frac{1}{[a_{n-1}]
 +\frac{1}{[a_n]}}} \]

 \noindent for $a_2, \ldots, a_n \in \ZZ - \{0\}$ and $n$ even or odd. We
 allow that the term
 $a_1$  may be zero, and in this case the tangle $[0]$ may be omitted.  A
 rational tangle described
 via a continued fraction in integer tangles is said to be in {\it
 continued
 fraction form}. The
 {\it length} of the continued fraction is arbitrary -- in the previous
formula illustrated
  with length $n$ -- whether the first summand is the tangle $[0]$ or not. }
 \end{defn}

 It follows  from Lemma 3.2 that inductively {\it  every rational
 tangle can be written in
 continued fraction form.}  Lemma 3.2 makes it easy to write out the
 continued
 fraction form of a given rational tangle, since horizontal twists are
 integer additions, and
 multiplications by vertical twists are the reciprocals of integer
 additions.
 For example, Figure 1
 illustrates the
  rational tangle $[2] + \frac{1}{[-2] +\frac{1}{[3]}},$  Figure 17
 illustrates the
  rational tangle $[2] + \frac{1}{[2] +\frac{1}{[3]}}.$  Note that

 \[ ([c] * \frac{1}{[b]}) + [a] \mbox{ \ has the continued fraction form \
 }
 [a]+ \frac{1}{[b]+
 \frac{1}{[c]}} = [[a], [b], [c]]. \]

 \noindent  For $T=[[a_1], [a_2],\ldots ,[a_n]]$ the following statements
 are
 now straightforward.

 \vspace{.1in}

 \noindent $\begin{array}{lrcl}

  1. & T + [\pm 1] &  = & [[a_1 \pm 1], [a_2],\ldots ,[a_n]],   \\
 [1.8mm]

  2. & \frac{1}{T} &  = & [[0], [a_1], [a_2],\ldots ,[a_n]],   \\  [1.8mm]

 3. & -T &  =  & [[-a_1], [-a_2],\ldots ,[-a_n]]. \\

 \end{array}$

\vspace{.1in}

  We now recall
 some facts about continued
 fractions. See for example
  \cite{Kh},  \cite{O}, \cite{Ko}, \cite{W}. In this paper we shall only
 consider continued fractions of the type

 \[ [a_1, a_2, \ldots, a_n] := a_1+ \frac{1}{a_2+\cdots + \frac{1}{a_{n-1}
 +\frac{1}{a_n}}}  \]

 \noindent for $ a_1 \in \ZZ, \ a_2, \ldots, a_n \in \ZZ - \{0\}$ and $n$
 even or odd.   The
 {\it length} of the continued fraction is the number
 $n$ whether $a_1$ is zero or not.   Note that if for $i>1$ all terms  are
 positive or all terms are negative and   $a_1 \neq 0$
 ($a_1 = 0,$) then the absolute value of the continued fraction is greater
 (smaller) than one. Clearly, the two simple algebraic operations {\it
addition of $+1$ or
 $-1$} and {\it inversion}
 generate inductively the whole class of continued fractions starting from
 zero.  For any rational number $\frac{p}{q}$ the following statements are
 straightforward.

 \vspace{.1in}

 \noindent \ $1.  \ \ \mbox{there are} \ a_1 \in \ZZ, \ a_2, \ldots, a_n
 \in
 \ZZ - \{0\} \
 \mbox{such that} \ \frac{p}{q} = [a_1, a_2, \ldots, a_n],  $

 \noindent $\begin{array}{lrcl}

  2.& \frac{p}{q} \pm 1 &  = & [a_1 \pm 1, a_2,\ldots ,a_n],   \\
 [1.8mm]

  3. & \frac{q}{p} &  = & [0, a_1, a_2, \ldots, a_n],   \\  [1.8mm]

 4. & -\frac{p}{q} &  =  & [-a_1, -a_2, \ldots, -a_n]. \\

 \end{array}$

\vspace{.1in}

 We can now define the fraction of a rational tangle.

 \begin{defn}{\rm \ Let T be a rational tangle isotopic to the continued fraction form
$[[a_1], [a_2],\ldots,[a_n]].$
 We define {\it the
 fraction $F(T)$ of $T$} to be the numerical  value of the continued
 fraction obtained by
 substituting integers for the  integer tangles in the expression for $T$,
 i.e.
 \[  F(T) := a_1+ \frac{1}{a_2+\cdots + \frac{1}{a_{n-1}
 +\frac{1}{a_n}}} = [a_1, a_2, \ldots, a_n],  \]
 \noindent if $T \neq [\infty],$ and $F([\infty]) := \infty = \frac{1}{0},$
 as a formal expression. 
  }
 \end{defn}

\begin{rem}{\rm \ This definition is good in the sense that one can show that  isotopic rational tangles always 
differ by flypes, and that the fraction is unchanged by flypes \cite{KL1}.   
} \end{rem}

 \noindent Clearly the tangle fraction has the following properties.

 \vspace{.1in}

 \noindent $\begin{array}{lrcl}

  1. & F(T + [\pm 1]) &  = & F(T) \pm 1,  \\
 [1.8mm]

  2. & F(\frac{1}{T}) &  = & \frac{1}{F(T)},   \\  [1.8mm]

 3. & F(-T) &  =  & -F(T). \\
 \end{array}$
 \vspace{.1in}

The main result about rational tangles (Theorem~1) is that two rational tangles are isotopic if and only if
they have the same fraction. We will show that every rational tangle is isotopic to a unique alternating 
continued fraction form, and that this alternating form can be deduced from the fraction of the tangle.
The Theorem then follows from this observation. 
\bigbreak

 \begin{lem}{ \  Every  rational tangle is isotopic to an alternating rational tangle.
 } \end{lem}

\noindent {\em Proof.}  Indeed, if $T$ has a non-alternating continued
fraction form then the following configuration, shown in the left of  Figure 11, must
occur somewhere in $T,$ corresponding to a change of sign from one term to an
adjacent term in the tangle continued fraction. This  configuration is isotopic to a
simpler isotopic configuration as shown in that figure.
\bigbreak

$$ \picill4.5inby1.05in(D11) $$
\begin{center}
 {\bf Figure 11 - Reducing to the alternating form}
\end{center}
\vspace{3mm}

 \noindent Therefore, it follows by induction on the number of crossings in
 the tangle that $T$ is isotopic to an alternating rational
 tangle.  $\hfill \Box$
 \bigbreak

Recall that a tangle is
alternating if and only if it has crossings all of
 the same type. Thus, {\it a rational tangle $T=[[a_1], [a_2],\ldots
,[a_n]]$
is  alternating if the $a_i$'s are all positive or all  negative}. For
example, the tangle of Figure 17 is alternating.
\bigbreak

A rational tangle $T=[[a_1], [a_2],\ldots ,[a_n]]$ is
 said to be in  {\it  canonical form} if $T$ is alternating and $n$ is odd.
The tangle of Figure 17 is in canonical form.
We  note that if $T$ is alternating  and $n$ even,
then we can bring $T$
 to canonical form by breaking $a_n$ by a unit, e.g.
 $[[a_1], [a_2], \ldots, [a_n]] = [[a_1], [a_2], \ldots, [a_n -1], [1]],$
 if $a_n >0.$  
\bigbreak

The last key observation is the following  well-known fact about continued
fractions.  

\begin{lem}{ \  Every  continued fraction $[a_1, a_2, \ldots, a_n]$ can be
 transformed to a unique canonical form $ [\beta_1, \beta_2, \ldots,
\beta_m],$ where
 all $\beta_i$'s are positive or all negative integers and $m$ is odd.
  } \end{lem}

 \noindent {\em Proof.}  It follows immediately from
 Euclid's algorithm. We evaluate first $[a_1, a_2, \ldots, a_n] =
 \frac{p}{q},$ and using
 Euclid's algorithm we rewrite $\frac{p}{q}$ in the desired form. We
 illustrate the proof with  an
 example. Suppose that
 $\frac{p}{q} =
 \frac{11}{7}$.  Then

 $$\frac{11}{7} = 1 + \frac{4}{7} = 1 + \frac{1}{\frac{7}{4}} = 1 +
 \frac{1}{1 + \frac{3}{4}} =
 1+ \frac{1}{1 + \frac{1}{\frac{4}{3}}} $$

 $$= 1+ \frac{1}{1 + \frac{1}{1 + \frac{1}{3}}} = [1,1,1,3] = 1+ \frac{1}{1
 +
 \frac{1}{1 + \frac{1}{2
 +
 \frac{1}{1}}}} = [1,1,1,2,1].$$

 \noindent This completes the proof. $\hfill \Box $
 \bigbreak

Note that if $T=[[a_1], [a_2],\ldots ,[a_n]]$ and $S = [[b_1], [b_2],\ldots ,[b_m]]$ are rational tangles in canonical form
with the same fraction, then it follows from this Lemma that $[a_1, a_2,\ldots ,a_n]$ and  $[b_1, b_2,\ldots ,b_m]$ are
canonical continued fraction forms for the same rational number, and hence are equal term-by-term. Thus the uniqueness of
canonical forms for continued fractions implies the uniqueness of canonical forms for rational tangles.  For
 example, let $T = [[2],[-3],[5]].$ Then $F(T) = [2, -3, 5] =
\frac{23}{14}.$
 But $ \frac{23}{14} =
  [1,1,1,1,4],$ thus $T \sim [[1],[1],[1],[1],[4]],$ and this last tangle
 is the canonical form of
 $T.$
\bigbreak

\noindent {\em Proof of Theorem~1.} We have now assembled all the ingredients for the proof of Theorem~1. In one direction,
suppose that rational tangles $T$ and $S$ are isotopic. Then each is isotopic to its canonical form $T'$ and $S'$ by a sequence of
flypes. Hence the alternating tangles $T'$ and $S'$ are isotopic to one another. By the Tait conjecture, there is a sequence of
flypes from 
$T'$ to $S'.$ Hence there is a sequence of flypes from $T$ to $S.$ One verifies that the fraction as we defined it is invariant
under flypes. Hence $T$ and $S$ have the same fraction. In the other direction, suppose that $T$ and $S$ have the same fraction.
Then, by the remark above, they have identical canonical forms to which they are isotopic, and therefore they are
isotopic to each other. This completes the proof of the Theorem. $\hfill \Box $

\section{Alternate Definitions of the Tangle Fraction} 

In the last section and in \cite{KL1} the
fraction of a rational tangle is defined directly from its combinatorial
structure, and we verify the
topological invariance of the fraction using the Tait conjecture.
\smallbreak

In \cite{KL1} we give yet another definition of the fraction for rational tangles
by using coloring of the tangle arcs.
There are definitions that associate a fraction
 $F(T)$ (including $0/1$ and $1/0$) to any $2$-tangle $T$ whether or not it
 is  rational. The first
 definition is due to John Conway in \cite{C1} using the Alexander
 polynomial of the knots $N(T)$ and
 $D(T).$ In \cite{GK2} an alternate definition is given that uses the
 bracket polynomial of the  knots
 $N(T)$ and $D(T),$ and in \cite{GK1}   the fraction of a tangle is related
 to the conductance of an
 associated electrical network.   
In all these
 definitions the fraction is by
 definition an isotopy invariant of tangles.   
Below we  discuss the bracket polynomial and coloring definitions of the fraction.

\subsection{$F(T)$ Through the Bracket Polynomial}

   In this section we shall discuss the structure of the the bracket state model for the Jones 
polynomial \cite{K1,K2} and how to construct the tangle fraction by using this technique.
We first construct the bracket polynomial (state summation), which is a regular isotopy invariant 
(invariance under all but the Reidemeister move I). The bracket polynomial can be normalized to produce an invariant of
all the Reidemeister moves. This invariant is known as the Jones polynomial \cite{Jones1,Jones2}. The Jones
polynomial was originally discovered by a different method. 
\bigbreak 

The {\em bracket polynomial} , $<K> \, = \, <K>(A)$,  assigns to each unoriented link diagram $K$ a 
Laurent polynomial in the variable $A$, such that
   
\begin{enumerate}
\item If $K$ and $K'$ are regularly isotopic diagrams, then  $<K> \, = \, <K'>$.
  
\item If  $K \amalg O$  denotes the disjoint union of $K$ with an extra unknotted and unlinked 
component $O$ (also called `loop' or `simple closed curve' or `Jordan curve'), then 

$$< K \amalg O> \, = \delta<K>,$$ 
where  $$\delta = -A^{2} - A^{-2}.$$
  
\item $<K>$ satisfies the following formulas 

$$<\mbox{\large $\chi$}> \, = A <\mbox{\large $\asymp$}> + A^{-1} <)(>$$
$$<\overline{\mbox{\large $\chi$}}> \, = A^{-1} <\mbox{\large $\asymp$}> + A <)(>,$$
\end{enumerate}

\noindent where the small diagrams represent parts of larger diagrams that are identical except  at
the site indicated in the bracket. We take the convention that the letter chi, \mbox{\large $\chi$},
denotes a crossing where {\em the curved line is crossing over the straight
segment}. The barred letter denotes the switch of this crossing, where {\em the curved
line is undercrossing the straight segment}. The above formulas can
be summarized by the single equation

$$<K> = A<S_{L}K> + A^{-1}<S_{R}K>.$$

\noindent  In this text formula we have used the notations $S_{L}K$ and $S_{R}K$ to indicate the 
two new diagrams created by the two smoothings of a single crossing in the diagram $K$. That is, $K$,
$S_{L}K$ and $S_{R}K$ differ at the site of one crossing in the diagram $K$.  These smoothings are
described as follows. Label the four regions locally incident to a crossing by the letters $L$ and
$R$, with $L$ labelling the region to the left of the undercrossing arc for a traveller who
approaches the overcrossing on a route along the undercrossing arc. There are two such routes, one
on each side of the overcrossing line. This labels two regions with $L$. The remaining two are
labelled $R$. A smoothing is of {\em type} $L$ if it connects the  regions labelled $L$, and it is
of {\em type} $R$ if it connects the regions labelled $R$, see Figure 12.

$$ \picill5inby2.8in(D12) $$
\begin{center} {\bf Figure 12 - Bracket Smoothings} 
\end{center}
\vspace{3mm}

\noindent It is easy to see that Properties $2$ and $3$ define the calculation of the bracket on
arbitrary link diagrams. The choices of coefficients ($A$ and $A^{-1}$) and the value of $\delta$
make the bracket invariant under the Reidemeister moves II and III (see \cite{K1}). Thus
Property $1$ is a consequence of the other two properties. 

\bigbreak

In order to obtain a closed formula for the bracket, we now describe it as a state summation.
Let $K$ be any unoriented link diagram. Define a {\em state}, $S$, of $K$  to be a choice of
smoothing for each  crossing of $K.$ There are two choices for smoothing a given  crossing, and
thus there are $2^{N}$ states of a diagram with $N$ crossings.
 In  a state we label each smoothing with $A$ or $A^{-1}$ according to the left-right convention 
discussed in Property $3$ (see Figure 12). The label is called a {\em vertex weight} of the state.
There are two evaluations related to a state. The first one is the product of the vertex weights,
denoted  

$$<K|S>.$$
The second evaluation is the number of loops in the state $S$, denoted  $$||S||.$$
  
\noindent Define the {\em state summation}, $<K>$, by the formula 

$$<K> \, = \sum_{S} <K|S>\delta^{||S||-1}.$$
It follows from this definition that $<K>$ satisfies the equations
  
$$<\mbox{\large $\chi$}> \, = A <\mbox{\large $\asymp$}> + A^{-1} <)(>,$$
$$<K \amalg  O> \, = \delta<K>,$$
$$<O> \, =1.$$
  
\noindent The first equation expresses the fact that the entire set of states of a given diagram is
the union, with respect to a given crossing, of those states with an $A$-type smoothing and those
 with an $A^{-1}$-type smoothing at that crossing. The second and the third equation
are clear from the formula defining the state summation. Hence this state summation produces the
bracket polynomial as we have described it at the beginning of the  section. 

\bigbreak

In computing the bracket, one finds the following behaviour under Reidemeister move I: 
  $$<\mbox{\large $\gamma$}> = -A^{3}<\smile> \hspace {.5in}$$ and 
  $$<\overline{\mbox{\large $\gamma$}}> = -A^{-3}<\smile> \hspace {.5in}$$

\noindent where \mbox{\large $\gamma$}  denotes a curl of positive type as indicated in Figure 13, 
and  $\overline{\mbox{\large $\gamma$}}$ indicates a curl of negative type, as also seen in this
figure. The type of a curl is the sign of the crossing when we orient it locally. Our convention of
signs is also given in Figure 13. Note that the type of a curl  does not depend on the orientation
we choose.  The small arcs on the right hand side of these formulas indicate
the removal of the curl from the corresponding diagram.  

\bigbreak
  
\noindent The bracket is invariant under regular isotopy and can be  normalized to an invariant of
ambient isotopy by the definition  
$$f_{K}(A) = (-A^{3})^{-w(K)}<K>(A),$$ where we chose an orientation for $K$, and where $w(K)$ is 
the sum of the crossing signs  of the oriented link $K$. $w(K)$ is called the {\em writhe} of $K$. 
The convention for crossing signs is shown in  Figure 13.

$$ \picill4.5inby2.5in(D13) $$
\begin{center} {\bf Figure 13 - Crossing Signs and Curls} 
\end{center}
\vspace{3mm}

\noindent By a change of variables one obtains the original
Jones polynomial, $V_{K}(t),$  for oriented knots and links from the normalized bracket:

$$V_{K}(t) = f_{K}(t^{-\frac{1}{4}}).$$

\noindent The bracket model for the Jones polynomial is quite useful both theoretically and in terms
 of practical computations. One of the neatest applications is to simply compute $f_{K}(A)$ for the
trefoil knot $T$ and determine that  $f_{T}(A)$ is not equal to $f_{T}(A^{-1}) = f_{-T}(A).$  This
shows that the trefoil is not ambient isotopic to its mirror image, a fact that is quite tricky to
prove by classical methods.
\bigbreak

For $2$-tangles, we do smoothings on the tangle diagram until there are no
crossings left. As a result, {\em a state of a $2$-tangle} consists in a collection of loops 
in the tangle box, plus simple arcs that connect the tangle ends. The loops evaluate to
powers of  $\delta$, and what is left is either the tangle $[0]$ or the tangle $[\infty]$, since  
$[0]$ and $[\infty]$ are the only ways to connect the tangle inputs and outputs without introducing
any crossings in the diagram.  In analogy to knots and links, we can find a {\em state
summation} formula for the {\em bracket of the tangle}, denoted $<T>,$ by summing over the states
obtained by smoothing each crossing in the tangle.  For this we define the {\em remainder of a
state}, denoted $R_S$, to be either the tangle $[0]$ or the tangle $[\infty]$. Then the evaluation
of $<T>$ is given by

$$<T> \, = \sum_{S} <T|S>\delta^{||S||}<R_S>,$$ where $<T|S>$ is the product of the vertex weights
($A$ or $A^{-1}$) of the state $S$ of $T$. The above formula is consistent with the
formula for knots obtained by taking  the closure $N(T)$ or $D(T)$. In fact, we have the following
formula:

$$<N(T)> \, = \sum_{S} <T|S>\delta^{||S||}<N(R_S)>.$$
Note that $<N([0])> = \delta$ and $<N([\infty])> = 1.$ A similar formula holds for $<D(T)>.$  Thus,
$<T>$ appears as a linear combination with Laurent polynomial coefficients of  $<[0]>$ and  
$<[\infty]>,$ i.e. $<T>$ takes values in the free module over $\ZZ[A,A^{-1}]$ with basis $\{ <[0]>,
 \, <[\infty]>\}.$ Notice that two elements in this module are equal iff the corresponding
coefficients of the basis elements coincide. Note also that $<T>$ is an invariant of regular isotopy
with values in this module. We have just proved the following:

\begin{lem}{ \  Let $T$ be any $2$-tangle and let $<T>$ be the formal expansion of the
bracket on this tangle. Then there exist elements $n_{T}(A)$ and $d_{T}(A)$ in $\ZZ[A,A^{-1}],$ such 
that

$$<T> \, = d_{T}(A)<[0]> + \, n_{T}(A)<[\infty]>,$$ and $n_{T}(A)$ and $d_{T}(A)$ are regular
isotopy invariants of the tangle $T$.
 } \end{lem}

\noindent In order to evaluate $<N(T)>$ in the formula above we need only apply the
closure $N$ to  $[0]$ and $[\infty].$ More precisely, we have: 

\begin{lem}{ \ $<N(T)> \, = d_T\delta + n_T$ and $<D(T)> \, = d_T + n_T\delta.$
 } \end{lem}

\noindent {\em Proof.} Since the smoothings of crossings do not interfere with the closure ($N$ or
$D$), the closure will carry through linearly to the whole sum of $<T>$. Thus, 

\vspace{.1in}

\noindent $\begin{array}{lcl}
\noindent <N(T)> & = & d_T(A)<N([0])> + n_T(A)<N([\infty])> = d_T(A) \delta + n_T(A),  \\   [1.8mm] 
\noindent <D(T)> & = &  d_T(A)<D([0])> + n_T(A)<D([\infty])> = d_T(A) + n_T(A) \delta. \hfill \Box 
 \end{array}$

\bigbreak

We define now the {\em polynomial fraction}, $frac_{T}(A)$, of the $2$-tangle $T$ to be the ratio

$$frac_{T}(A) = \frac{n_T(A)}{d_T(A)}$$

\noindent in the ring of fractions of $\ZZ[A,A^{-1}]$ with a formal symbol $\infty$ adjoined.

\begin{lem}{ \  $frac_{T}(A)$ is an invariant of ambient isotopy for $2$-tangles.
 } \end{lem}

\noindent {\em Proof.} Since $d_{T}$ and $n_{T}$ are regular isotopy invariants of $T$,  it
follows that  $frac_{T}(A)$ is also a regular isotopy invariant of $T$. Suppose now $T\gamma$ is
$T$ with  a curl added. Then $<T\gamma> \, = (-A^3)<T>$ (same remark for $\bar{\gamma}$). So,
$n_{T\gamma}(A)  = -A^3 n_T(A)$ and $d_{T\gamma}(A) = -A^3 d_T(A).$ Thus, 
${n_{T\gamma}}/{d_{T\gamma}} = {n_{T}}/{d_{T}}.$ This shows that  $frac_{T}$ is also
invariant under the Reidemeister move I, and hence an  ambient isotopy invariant. $\hfill \Box$

\begin{lem}{ \  Let $T$ and $S$ be two $2$-tangles. Then, we have the following formula for the
bracket of the sum of the tangles.

$$<T+S> \, = d_{T}d_{S}<[0]> + \, (d_{T}n_{S} + n_{T}d_{S} + n_{S}\delta)<[\infty]>.$$

\noindent Thus $$frac_{T+S} = frac_{T} + frac_{S} + \frac{n_{S}\delta}{d_{T}d_{S}}.$$
 } \end{lem}

\noindent {\em Proof.} We do first the smoothings in $T$ leaving $S$ intact, and then in $S$: 

\vspace{.1in}

\noindent $\begin{array}{rcl}
<T+S> & = & d_{T}<[0] + S> + n_{T}<[\infty] + S>    \\   [1.8mm] 
   \ &  = &  d_{T}<S> +  n_{T}<[\infty] + S>        \\   [1.8mm] 
   \ &  = &  d_{T}(d_{S}<[0]> + n_{S}<[\infty]>)    \\   [1.8mm] 
   \ &  + &  n_{T}(d_{S}<[\infty] + [0]> + n_{S}<[\infty] + [\infty]>)    \\   [1.8mm] 
   \ &  = &  d_{T}(d_{S}<[0]> + n_{S}<[\infty]>) + n_{T}(d_{S}<[\infty]> + n_{S} \delta<[\infty]>)
                \\   [1.8mm] 
   \ &  = &   d_{T}d_{S}<[0]> + (d_{T}n_{S} + n_{T}d_{S} + n_{S} \delta) <[\infty]>.   \\     
 \end{array}$

\vspace{.1in}
\noindent Thus, $n_{T+S} = (d_{T}n_{S} + n_{T}d_{S} + n_{S} \delta)$ and 
$d_{T+S} = d_{T}d_{S}.$ A straightforward calculation gives now $frac_{T+S}.  \hfill \Box $

\bigbreak
\noindent As we see from Lemma 4, $frac_{T}(A)$ will be additive on tangles if

$$\delta = -A^{2} - A^{-2} =0.$$

\noindent Moreover, from Lemma 2 we have  for $\delta = 0, \ <N(T)> \, = n_T, \ <D(T)> \, = d_T.$
 This nice situation will be our main object of study in the rest of this section. Now, if we set 
 $A=\sqrt{i}$ where $i^{2}=-1,$ then it is 

$$\delta = -A^{2} -A^{-2} = -i - i^{-1} = -i +i = 0.$$

\noindent For this reason, we shall henceforth assume that $A$ takes the value $\sqrt{i}$. So $<K>$
will denote $<K>(\sqrt{i})$ for any knot or link $K$. 

\bigbreak

We now define the {\em $2$-tangle fraction} $F(T)$ by the following formula:

$$F(T) = i\, \frac{n_{T}(\sqrt{i})}{d_{T}(\sqrt{i})}.$$

\noindent We will let $n(T) = n_{T}(\sqrt{i})$ and $d(T) = d_{T}(\sqrt{i}),$ so that 

$$F(T) = i\, \frac{n(T)}{d(T)}.$$

\begin{lem}{ \  The $2$-tangle fraction has the following properties.
\begin{enumerate}
\item  \ \ \  $F(T) = i\, {<N(T)>}/{<D(T)>},$ and it is a real number or $\infty$,

\item  \ \ \  $F(T+S) = F(T) + F(S),$

\item  \ \ \  $F([0]) = \frac{0}{1},$

\item  \ \ \  $F([1]) = \frac{1}{1},$

\item  \ \ \  $F([\infty]) = \frac{1}{0},$

\item  \ \ \  $F(-T) = - F(T), \ \mbox{in particular} \ F([-1]) = -\frac{1}{1},$

\item  \ \ \  $F({1}/{T}) = {1}/{F(T)},$

\item  \ \ \  $F(T^r) = {-1}/{F(T)}.$
\end{enumerate}

 } \end{lem}

\noindent As a result we conclude that for a tangle obtained by arithmetic operations from
integer tangles $[n],$ the fraction of that tangle is the same as the
fraction obtained by doing the same operations to the corresponding integers. (This will be studied
in detail in the next section.) 

\bigbreak
\noindent
{\em Proof.} The formula $F(T) = i\, {<N(T)>}/{<D(T)>}$ and Statement $2.$ follow from the
observations above about $\delta = 0.$ In order to show that $F(T)$ is a real number or $\infty$ we
first consider $<K> := <K>(\sqrt{i}),$ for $K$ a knot or link, as in the hypotheses prior to the
lemma. Then we apply this information to the ratio $i\, {<N(T)>}/{<D(T)>}.$ 

Let $K$ be any knot or link.  We claim that then $<K> = \omega p,$ where $\omega$ is a power of 
$\sqrt{i}$ and $p$ is an integer. In fact, we will show that each non-trivial state of $K$ 
contributes $\pm \omega$ to $<K>.$ In order to show this, we examine how to get from one
non-trivial state to another. It is a fact that, for any two states, we can get from one to the
other by resmoothing a subset of crossings. It is possible to get from any single loop state (and
only single loop states of $K$ contribute to  $<K>$, since $\delta = 0$) to any other single loop
state by a series of {\it double resmoothings}. In a double resmoothing we resmooth two crossings,
such that one of the resmoothings disconnects the state and the other reconnects it. See Figure 14
for an illustration. Now consider the effect  of a double resmoothing on the 
evaluation of one state. Two crossings change. If one is labelled $A$ and the other $A^{-1}$, then 
there is no net change in the  evaluation of the state. If both are $A$, then we go from $A^{2}P$
($P$  is the rest of the product of state labels) to $A^{-2}P.$ But $A^2 = i$ and  $A^{-2} = -i.$
Thus if one state contributes $\omega = ip$, then the other state contributes $-\omega = -ip.$ These
remarks prove the claim. 

$$ \picill2inby2in(D14) $$
\begin{center} {\bf Figure 14 - A Double Resmoothing} 
\end{center}
\vspace{3mm}

Now, a state that contributes non-trivially to $N(T)$ must have the form of the tangle $[\infty].$
We will show that if $S$ is a state of $T$ contributing non-trivially to $<N(T)>$ and $S'$ a state
of $T$ contributing non-trivially to $<D(T)>,$ then ${<S>}/{<S'>} = \pm i.$ Here $<S>$ denotes the
product of the vertex weights for $S$, and $<S'>$ is the product of the vertex weights for $S'.$ If
this ratio is verified for some pair of states $S,\, S',$ then it follows from the first claim
that it is true for all pairs of states, and that $<N(T)> = \omega p,$ \ $<D(T)> = {\omega}' q, \
p, q \in \ZZ$ and $\omega/{\omega}' = {<S>}/{<S'>} = \pm i.$ Hence  ${<N(T)>}/{<D(T)>} = \pm i\,
p/q,$ where $p/q$ is a rational number (or $q=0$). This will complete the proof that $F(T)$ is real
or $\infty.$

To see this second claim we consider specific pairs of states as in Figure 15. We have
illustrated representative states  $S$ and $S'$ of the tangle $T$. We obtain  $S'$ from $S$ by
resmoothing at one site that changes $S$ from an $[\infty]$ tangle to the $[0]$ tangle underlying
$S'$. Then  ${<S>}/{<S'>} = A^{\pm 2} = \pm i.$ If there is no such resmoothing site available, then
it follows that $D(T)$ is a disjoint union of two diagrams, and hence $<D(T)> \, =0$ and $F(T)
=\infty.$ This does complete the proof of Statement 1. 

$$ \picill3inby3in(D15) $$
\begin{center} {\bf Figure 15 - Non-trivial States} 
\end{center}
\vspace{3mm}

At $\delta = 0$ we also have: 

\vspace{.1in}

\noindent $<N([0]> \, =0, \, <D([0])> \, = 1, \, <N([\infty])> \, = 1, \, <D([\infty])> = \, 0,$ and
so, the evaluations $3.$ to $5.$ are easy. For example, note that

$$<[1]> \, = A<[0]> + A^{-1}<[\infty]>,$$

\noindent hence $$F([1]) = i\, \frac{A^{-1}}{A} = i\, A^{-2} = i\, (i^{-1}) = 1.$$

\noindent To have the fraction value $1$ for the tangle $[1]$ is the reason that in the definition
of $F(T)$ we normalized by $i$. Statement $6.$ follows from the fact that the bracket of the mirror
image of a knot $K$ is the same as the bracket of $K$, but with $A$ and $A^{-1}$ switched.
For proving $7.$ we observe first that for any $2$-tangle
$T$, $d(\frac{1}{T}) = \overline{n(T)}$ and $n(\frac{1}{T})= \overline{d(T)},$ where the
overline denotes the complex conjugate. Complex conjugates occur because $A^{-1} = \overline{A}$
when $A = \sqrt{i}$. Now, since $F(T)$ is real, we have 

\vspace{.1in}
\noindent $F(\frac{1}{T}) = i \, {\overline{d(T)}}/{\overline{n(T)}} = 
\overline{-i \, {d(T)}/{n(T)}} =\overline{{1}/({i \, {n(T)}/{d(T)}})} =
\overline{{1}/{F(T)}} = {1}/{F(T)}.$
\vspace{.1in}

\noindent Statement $8.$ follows immediately from $6.$ and $7.$ This completes the proof. $\hfill
\Box $

\bigbreak

For a related approach to the well-definedness of the $2$-tangle fraction, the reader should consult \cite{Krebes}.
The double resmoothing idea originates from \cite{K5}. 

\begin{rem}{\rm \ For any knot or link $K$ we define the {\em determinant} of $K$ by the 
formula

$$Det(K) := |<K>(\sqrt{i})|$$

\noindent where $|z|$ denotes the modulus of the complex number $z$. Thus we have the formula

$$|F(T)| = \frac{Det(N(T))}{Det(D(T))}$$

\noindent for any $2$-tangle $T$.  

In other approaches to the theory of knots, the determinant of 
the knot is actually the determinant of a certain matrix associated either to the diagram for the
knot or to a surface whose boundary is the knot. See \cite{R, K4} for more information
on these connections. Conway's original definition of the fraction \cite{C1} is $\Delta_{N(T)}(-1)/\Delta_{D(T)}(-1)$
where $\Delta_{K}(-1)$ denotes the evaluation of the Alexander polynomial of a knot $K$ at the value $-1.$ In fact,
$|\Delta_{K}(-1)| = Det(K),$ and with appropriate attention to signs, the Conway definition and our definition using the
bracket polynomial coincide for all $2$-tangles.   
 } \end{rem}

\subsection{The Fraction through Coloring}
 
We conclude this section by giving an alternate definition of the fraction
 that uses the concept of
 coloring of knots and tangles. We color the arcs  of the knot/tangle with
 integers, using the basic
 coloring rule  that if two undercrossing arcs colored $\alpha$ and
 $\gamma$  meet at an overcrossing arc colored $\beta$, then \, $\alpha +
 \gamma = 2\beta.$
 We often think of one of the undercrossing arc colors as determined
 by  the other two  colors. Then one writes $\gamma = 2\beta -\alpha.$
 \smallbreak

 It is easy to verify
 that this coloring method  is invariant under  the Reidemeister moves in
 the following sense: Given a
 choice of coloring for the tangle/knot, there is a way to re-color it
 each time a Reidemeister
 move is performed, so that no change occurs to the colors on the external
 strands of the tangle (so
 that we still have a valid coloring). This means that a  coloring
 potentially contains topological
 information about a knot or a tangle. In coloring
 a knot (and also many non-rational tangles) it is usually necessary to
 restrict the colors to the set
 of integers modulo $N$ for some modulus $N$. For example, in  Figure 16 it
 is clear
 that the color set ${\ZZ}/{3\ZZ} =
 \{ 0,1,2 \}$ is forced for coloring a trefoil knot.  When there exists a
 coloring of a tangle by
 integers, so that it is not necessary to reduce the colors over some
 modulus we shall say that the tangle is {\it integrally colorable}.
\bigbreak

$$ \picill4.2inby1.18in(D16) $$
 \begin{center}
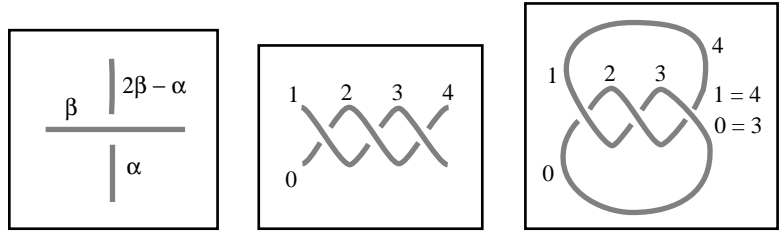
{\bf Figure 16 - The coloring rule, integral and modular coloring }
 \end{center}
 \vspace{3mm}

It turns out that {\it every rational tangle is integrally colorable:} To see this
 choose two `colors' for the initial strands (e.g. the colors $0$ and $1$)
and color
 the rational tangle as you
 create  it by successive twisting. We call the colors on the initial
 strands the {\it starting
 colors}.  See Figure 17 for an example. It is important that we start
coloring from the initial
strands, because then the coloring propagates automatically and uniquely. If
one starts from
 somewhere else, one might get into an edge with an undetermined color. The
 resulting colored tangle now has colors assigned to its external strands at
the northwest,
 northeast, southwest and southeast positions.  Let
 $NW(T)$, $NE(T)$, $SW(T)$ and
 $SE(T)$ denote these respective colors of the colored tangle $T$ and define
 the {\em color matrix of $T$}, $M(T)$, by the equation

 $$M(T) =  \left[
 \begin{array}{cc}
      NW(T)& NE(T) \\
      SW(T) & SE(T)
 \end{array}
 \right]. $$

 \begin{defn}{\rm \  To  a rational tangle $T$  with color matrix
 $M(T) =  \left[
 \begin{array}{cc}
      a& b \\
      c & d
 \end{array}
 \right]$   we associate the number
 \[ f(T) := \frac{b-a}{b-d} \ \in \QQ \cup \infty. \]
 }
 \end{defn}

 \noindent It turns out that the entries $a, b, c, d$ of a color matrix of
 a
 rational tangle  satisfy
 the   `diagonal sum rule': \ $a+d=b+c$.

\bigbreak

$$ \picill4.3inby2.25in(D17) $$
\begin{center}
{\bf Figure 17 - Coloring Rational Tangles }
 \end{center}
 \vspace{3mm}

 \begin{prop}{ \ The number $f(T)$ is a topological invariant associated
 with

 the tangle $T$.  In fact, $f(T)$ has  the following properties:
 \vspace{.1in}

 \noindent $\begin{array}{lrcl}

  1. & f(T+[\pm1]) &  = & f(T) \pm 1,  \\    [1.8mm]

  2. & f(-\frac{1}{T})  &  = & -\frac{1}{f(T)},   \\  [1.8mm]

 3. & f(-T) &  =  & -f(T),  \\  [1.8mm]

 4. & f(\frac{1}{T}) &  =  & \frac{1}{f(T)},  \\  [1.8mm]

 5. & f(T) &  =  & F(T).  \\
 \end{array}$
\vspace{.1in}

 \noindent Thus the coloring fraction is identical to the arithmetical
 fraction defined earlier.
  } \end{prop}

 \noindent  It is easy to see that $f([0])
 = \frac{0}{1}$, $ f([\infty])= \frac{1}{0}$, $f([\pm 1])=
 \pm 1.$ Hence Statement 5 follows by  induction. For proofs of all
statements above
 as well as for a more
 general set-up we refer the reader to our paper
 \cite{KL1}.  This  definition is quite elementary, but applies only to
 rational
 tangles and tangles generated from them by the algebraic operations of
 `$+$'
 and `$*$'.
\bigbreak

In Figure 17 we have illustrated a coloring over the integers for the tangle
$[[2],[2],[3]]$ such that every edge is labelled by a different integer. 
This is always the case for an alternating rational tangle diagram $T.$
For the numerator closure $N(T)$ one obtains a coloring in a modular number system.
For example in Figure~17 the coloring of $N(T)$ will be in $\ZZ/17\ZZ$, and it is easy
to  check that the labels remain distinct in this example. 
For rational tangles, this is always the case when $N(T)$ has a prime 
determinant, see \cite{KL1} and \cite{LP}. It is part of a more general conjecture about alternating knots and links
\cite{KH, APA}.

\subsection{The Fraction through Conductance} 

Conductance is a quantity defined in electrical networks as the inverse of resistance. For pure resistances, conductance is a
positive quantity.  Negative conductance
corresponds to amplification, and is commonly included in the physical formalism. One defines the conductance
between two vertices in a graph (with positive or negative conductance weights on the  edges of the graph) as a sum of weighted
trees in the graph divided by a sum of weighted trees of the same graph, but with the two vertices identified. This definition
allows negative values for conductance and it agrees with the classical one.   Conductance satisfies familiar laws of  parallel
and series connection as well as a star-triangle relation. 
\bigbreak

By associating to a given knot or link diagram the corresponding signed checkerboard
graph (see \cite{KL1, GK1} for a definition of this well-known association of graph to link diagram), one can define \cite{GK1} the
conductance of a knot or link between any two regions that receive the same color in the checkerboard graph. The conductance of
the link between these two regions is an isotopy invariant of the link (with motion restricted to Reidemeister moves that do not
pass across the selected regions). This invariance follows from properties of series/parallel connection and the star-triangle
relation. These circuit laws turn out to be images of the Reidemeister moves under the translation from knot or link diagram to
checkerboard graph! For a $2$-tangle we take the conductance to be the conductance of the numerator of the tangle, between the
two bounded regions adjacent to the closures at the top and bottom of the tangle.   
\bigbreak

The conductance of a $2$-tangle turns out to be the same as the fraction of the tangle. This provides yet another
way to define and verify the isotopy invariance of the tangle fraction for any $2$-tangle.

 \section{The Classification of Unoriented Rational Knots}

 By taking their numerators or denominators rational tangles give rise to a
 special class of knots,
 the rational knots.  We have seen so far that rational tangles
 are directly related to finite continued fractions. We carry
 this
 insight further into the
 classification of rational knots (Schubert's theorems). In this section we
 consider unoriented knots,
 and by Remark 3.1 we will be  using the
 $3$-strand-braid representation for rational tangles with odd number of
 terms. Also, by Lemma 2, we may assume all rational knots to
 be alternating.
 Note that we  only need
 to take numerator closures, since the denominator closure of a tangle is
 simply the numerator closure of its rotate.
 \smallbreak

 As already said in the introduction, it may happen that two
 rational tangles are non-isotopic but have
 isotopic  numerators. The simplest instance of this phenomenon is adding
 $n$ twists at the
 bottom of a tangle $T$, see Figure 18. This operation does not change the
 knot $N(T),$ i.e.
 $N(T*1/[n]) \sim N(T),$ but it does change the tangle, since 
$F(T*1/[n]) = F(1/([n] + 1/T)) = 1/(n + 1/F(T));$ so, if $F(T) = p/q$, then
 $F(T*1/[n]) = p/(np + q).$ Hence, if we set $np + q = q'$ we have $q\equiv
 q'(mod \, p),$ just as Theorem~2 dictates. Note that
reducing all possible bottom twists implies $|p|>|q|.$
\bigbreak

 $$ \picill3inby1.6in(D18) $$
\begin{center}
{\bf Figure 18 - Twisting the bottom of a tangle}
\end{center}
\vspace{3mm}

 Another key example of the arithmetic relationship of the
 classification of rational knots is illustrated in  Figure 19. Here we see
 that the `palindromic'
 tangles

 $$ T=[[2],[3],[4]] = [2] + \frac{1}{[3] + \frac{1}{[4]}} $$

 \noindent and
 $$ S=[[4],[3],[2]] = [4] + \frac{1}{[3] + \frac{1}{[2]}} $$

 \noindent both close to the same rational knot, shown at the bottom of the
 figure.
 The two tangles are different, since they have different corresponding
 fractions:

  $$ F(T)=2 + \frac{1}{3 + \frac{1}{4}} = \frac{30}{13}  \mbox{ \ \ and \ \
 }
 F(S)=4 + \frac{1}{3 +
 \frac{1}{2}} = \frac{30}{7}.$$ Note that the product of $7$ and $13$ is
 congruent to $1$ modulo $30.$
\bigbreak

 $$ \picill4.5inby2.45in(D19) $$
\begin{center}
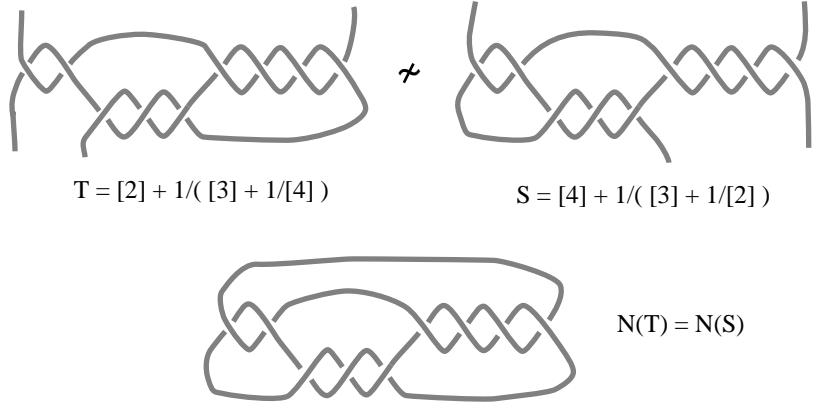
{\bf Figure 19 - An Instance of the Palindrome Equivalence}
 \end{center}
 \vspace{3mm}

 \noindent   More generally, consider the
 following two fractions:

 $$ F = [a,b,c] = a+ \frac{1}{b+ \frac{1}{c}}  \mbox{ \ \ and \ \ } G =
 [c,b,a] = c+ \frac{1}{b+
 \frac{1}{a}}.$$

 \noindent  We find that

 $$F = a + c\, \frac{1}{cb+1}  = \frac{abc+a +c}{bc+1} = \frac{P}{Q},$$

 \noindent while

 $$G = c + a\, \frac{1}{ab +1} = \frac{abc+c+a}{ab+1} = \frac{P}{Q'}.$$

 \noindent Thus we found that  $F= \frac{P}{Q}$ and $G= \frac{P}{Q'},$
 where

 $$QQ' = (bc+1)(ab+1) = ab^{2}c +ab +bc +1  = bP+1.$$

 \noindent Assuming that $a$, $b$ and $c$ are integers, we conclude that

 $$QQ'  \equiv  1 \, (mod\, P).$$

 \noindent  This pattern generalizes to arbitrary continued fractions and
 their palindromes (obtained by reversing the  order of the terms).
 I.e. {\it If $\{ a_{1},a_{2},\ldots ,a_{n} \}$ is a collection of $n$
 non-zero integers, and if $A = [a_{1},a_{2},\ldots ,a_{n}] = \frac{P}{Q}$
and
 $B~=~[a_{n},a_{n-1},\ldots ,a_{1}]= \frac{P'}{Q'},$ then  $P=P'$ and  $QQ'
 \equiv
 (-1)^{n+1} (mod\, P).$}  We will be referring to this as `the Palindrome
 Theorem'. The Palindrome Theorem is a known result about continued
fractions. For example, see \cite{Sie} and \cite{KL2}. 
  Note  that we need $n$ to be odd in the previous congruence.
 This agrees with Remark 3.1 that without loss of generality the terms in the
 continued fraction of a rational tangle may be assumed to be odd.
 \bigbreak

 Finally, Figure 20 illustrates another basic example for the
 unoriented Schubert Theorem. The two tangles
 $R = [1] + \frac{1}{[2]} $ and $S = [-3]$ are non-isotopic by the Conway
 Theorem, since $F(R) = 1 +
 1/2 = 3/2$ while $F(S) = -3 = 3/-1.$ But they have isotopic
 numerators:  $N(R) \sim N(S),$ the
 left-handed trefoil. Now  $2$ is congruent to $-1$ modulo $3,$  confirming
 Theorem~2.
\bigbreak

$$ \picill3.3inby1.45in(D20) $$
\begin{center}
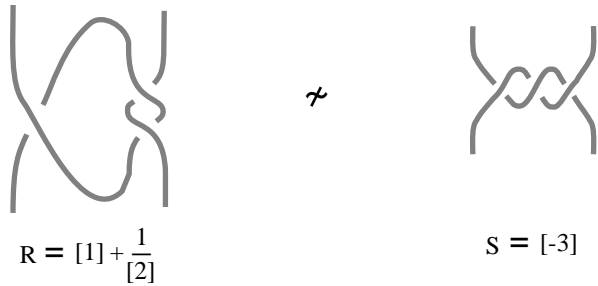
{\bf Figure 20 - An Example of the Special Cut}
 \end{center}
 \vspace{3mm}

We now analyse the above example in general. From the analysis
of the bottom twists we can assume without loss of generality that
a rational
tangle
$R$ has fraction
$\frac{P}{Q},$ for
 $|P|>|Q|.$ Thus $R$ can be
 written in the form $R=[1]+T$ or $R=[-1]+T.$ We consider the rational knot
 diagram $K = N([1]+T),$
 see Figure 21. (We analyze $N([-1]+T)$ in the same way.)  The tangle
 $[1]+T$ is said to arise as a {\it standard cut} on $K.$
\bigbreak

 $$ \picill5.2inby1.05in(D21) $$
\begin{center}
{\bf Figure 21 - A Standard Cut }
 \end{center}
 \vspace{3mm}

 Notice that the indicated horizontal crossing of $N([1]+T)$
 could be also seen as a
 vertical one. So, we could also cut the diagram $K$ at the two other
 marked points (see Figure 22) and still obtain
 a rational tangle, since $T$ is rational.  The tangle
 obtained by cutting $K$ in this
 second pair of points is said to arise as a {\it special cut} on
 $K.$  Figure 22 demonstrates that the tangle of the special cut is the
 tangle $[-1] - 1/T.$ \, So we have
 $N([1]+T) \sim N([-1] -\frac{1}{T}).$ Suppose  now $F(T) = p/q.$ Then
 $F([1]+T) = 1+ p/q = (p+q)/q,$
 while $F([-1] - 1/T) = -1 - q/p = (p+q)/(-p),$ so the two rational tangles
 that give rise to the same
 knot $K$ are not isotopic. Since $-p \equiv q\, mod  (p+q),$ this
 equivalence
 is another example for
 Theorem~2.  In Figure 22 if we took $T= \frac{1}{[2]} $ then
 $[-1] - 1/T = [-3]$ and we would obtain the example of Figure 20.
\bigbreak

$$ \picill5.35inby2.7in(D22) $$
\begin{center}
{\bf Figure 22 - A Special Cut }
 \end{center}
\vspace{3mm}

 The proof of  Theorem~2 can now proceed in two
 stages.   First, given a rational knot diagram we look for all possible
 places where we could cut
 and open it to a rational tangle.  The crux of our  proof in
 \cite{KL2} is the fact that  all possible `rational cuts' on a rational
 knot
 fall into one of the
 basic cases that we have already discussed. I.e. we have the  standard
 cuts,
 the palindrome cuts and
 the special cuts. In Figure 23 we
 illustrate on a representative rational knot, all the cuts that exhibit
 that
 knot as a closure of a
 rational tangle. Each pair of points is marked with the same number.
  The arithmetics is similar to the cases that have been already  verified.
 It is
 convenient to  say that reduced fractions
 $p/q$ and
 $p'/q'$ are {\it arithmetically equivalent}, written  $p/q \sim p'/q'$ if
 $p
 = p'$  and either $qq'
 \equiv 1$ (mod $p$)~or
~$q \equiv q'$~(mod~$p$~)~. In this language, Schubert's theorem  states that
two
 rational tangles close to form
 isotopic knots if and only if their fractions are arithmetically
 equivalent.
\bigbreak

 $$ \picill4.7inby2.5in(D23) $$
\begin{center}
{\bf Figure 23 - Standard, Palindrome and Special Cuts}
\end{center}
 \vspace{3mm}

 \noindent In Figure 24 we illustrate one example of a cut that is not
 allowed  since it opens the
 knot to a non-rational tangle.
\bigbreak

$$ \picill4.65inby1.45in(D24) $$
 \begin{center}
{\bf Figure 24 - A Non-Rational Cut}
 \end{center}
\vspace{3mm}

 In the second stage of the proof  we want to check the
 arithmetic equivalence for  two
 different given knot diagrams, numerators of some rational tangles. 
 By Lemma 2 the two knot
 diagrams may be assumed alternating, so by the Tait Conjecture they will
 differ by flypes. We analyse
 all possible flypes to prove that no new cases for study arise.  Hence the
 proof becomes complete at that point.  We
 refer the reader to our paper \cite{KL2} for the details.  $\hfill \Box $

\begin{rem} \rm \ The original proof of the classification of unoriented rational knots by Schubert \cite{Sch}
proceeded by a different route than the proof we have just sketched. Schubert used a $2$-bridge representation of 
rational knots (representing the knots and links as diagrams in the plane with two special overcrossing arcs, called the
bridges). From the $2$-bridge representation, one could extract a fraction $p/q,$ and Schubert showed by means of a canonical form,
that if two such  presentations are isotopic, then their fractions are arithmetically equivalent (in the sense that we have
described here). On the other hand, Seifert \cite{Sch} observed that the $2$-fold branched covering space of a $2$-bridge
presentation with fraction
$p/q$ is a lens space of type $L(p,q).$ Lens spaces are a particularly tractable set of three manifolds, and it is known by work
of  Reidemeister and Franz \cite{Rd3, Franz} that $L(p,q)$ is homeomorphic to $L(p',q')$ if and only if $p/q$ and $p'/q'$ are
arithmetically  equivalent. Furthermore, one knows that if knots $K$ and $K'$ are isotopic, then their $2$-fold branched
covering spaces are homeomorphic. Hence it follows that if two rational knots are isotopic, then their fractions are
arithmetically equivalent (via the result of Reidemeister and Franz classifying lens spaces). In this way Schubert proved that two
rational knots are isotopic if and only if their fractions are arithmetically equivalent.
\end{rem}

 \section{Rational Knots and Their Mirror Images}

 In this section we give an application of Theorem~2.
 An unoriented knot or link $K$ is said to be {\em achiral} if it is
 topologically equivalent to its
 mirror image $-K$. If a link is not equivalent to its mirror image then it
 is said be {\em chiral}.
 One then can speak of the  {\em chirality} of a given knot or link,
 meaning whether it is chiral or
 achiral. Chirality plays an important role in the applications of Knot
 Theory to Chemistry and
 Molecular Biology. It is interesting to use the classification of
 rational knots and links to
 determine their chirality. Indeed, we have the following well-known result
(for example see \cite{Sie} and also
 page 24, Exercise 2.1.4 in \cite{Kaw}):

 \begin{th} {\ Let $K=N(T)$ be an unoriented rational knot or link,
 presented as the numerator of a
 rational  tangle $T.$  Suppose that  $F(T) = p/q$ with $p$ and $q$
 relatively prime. Then $K$  is achiral if and only if
  $q^{2}\equiv -1 \, (mod \, p).$  It follows that achiral rational knots
 and links are all numerators of rational tangles  of the form
 $[[a_1], [a_2], \ldots, [a_k], [a_k], \ldots, [a_2], [a_1]]$ for any
 integers $a_1, \ldots, a_k.$ }
 \end{th}

\noindent  Note that in this description we are using
 a representation of the tangle with an even number of terms. The leftmost
twists $[a_1]$
 are horizontal, thus the rightmost starting twists  $[a_1]$ are vertical.
\bigbreak

 \noindent {\em Proof.}   With $-T$ the mirror image of the tangle $T$, we
 have that $-K = N(-T)$ and
 $F(-T) = p/(-q).$ If $K$ is topologically equivalent to $-K$, then $N(T)$
 and $N(-T)$ are
 equivalent, and it follows from the  classification theorem for rational
 knots that either $q(-q)
 \equiv 1 \, (mod \, p)$ or $q \equiv -q \, (mod \, p).$ Without loss of
 generality we can assume
 that $0< q < p.$ Hence $2q$ is not divisible by $p$ and therefore it is
 not the case that  $q \equiv
 -q \, (mod \, p).$ Hence $q^{2} \equiv -1 \, (mod \, p).$
 \smallbreak
 \noindent Conversely, if $q^{2} \equiv -1 \,
 (mod \, p),$ then it follows from the Palindrome Theorem (described in the previous section) \cite{KL2} that {\it the
 continued fraction expansion of
 $p/q$ has to be symmetric with an even number of terms.} It is then easy
 to see that the
 corresponding rational knot or link, say $K = N(T),$  is equivalent to its
mirror image. One
 rotates $K$ by $180^{\circ}$ in the plane and swings an arc, as Figure 25
illustrates. This completes the proof.
 $\hfill \Box$
\bigbreak

In \cite{ES0} the authors find an explicit formula for the number of achiral rational
knots among all rational knots with $n$ crossings.

$$ \picill4.8inby2.05in(D25) $$
\begin{center}
 {\bf Figure 25 - An Achiral Rational Link}
 \end{center}
\bigbreak

 \section{The Oriented Case}

 Oriented rational knots and links arise as numerator closures of oriented
 rational tangles.  In  order
 to compare oriented rational knots via rational tangles we need to examine
 how rational tangles
 can  be oriented. We orient rational tangles by choosing an orientation
 for each strand of the tangle. Here we are only interested in orientations
 that yield
 consistently oriented knots  upon taking the numerator closure. This means
 that the two top end arcs
 have to be oriented one inward and the other outward. Same for the two
 bottom end arcs.
 We shall say that two oriented rational tangles are {\it isotopic} if they
 are isotopic as unoriented
 tangles, by an isotopy that carries the orientation of one tangle to the
orientation of the
 other.  Note that, since the end arcs of a tangle are fixed during a tangle
isotopy, this means that
 the tangles must have identical orientations at their four end arcs {\it
NW, NE, SW, SE}.
 It follows that if we change the orientation of one or both strands of an
oriented rational tangle
 we will always obtain a non-isotopic oriented rational tangle.
\bigbreak

Reversing the orientation of one strand of an oriented rational
tangle may or may not give
rise to isotopic oriented rational knots. Figure 26 illustrates an example
of non-isotopic oriented
rational knots, which are isotopic as unoriented knots.
\bigbreak

 $$ \picill2inby2in(D26) $$
 \begin{center}
 {\bf Figure 26 -  Non-isotopic Oriented Rational Links }
 \end{center}
 \vspace{3mm}

 Reversing the orientation of both strands of an oriented
 rational tangle will always give rise
 to two isotopic oriented rational knots or links. We can see this by doing
a vertical flip, as Figure
 27 demonstrates.  Using this observation we conclude that, as far as the
study of oriented
 rational knots is concerned, {\it all oriented rational
 tangles  may be assumed to have the same orientation for their  NW and
 NE end arcs.}  We fix this orientation to be downward for the  {\it NW} end
arc and upward
 for the  {\it NE} arc, as in the examples of Figure 26 and as illustrated
in Figure 28.
 Indeed, if the orientations are
 opposite of the fixed ones doing
 a vertical flip  the knot may be considered as the
 numerator of the vertical flip
 of the original tangle. But this is unoriented isotopic to the original
tangle
 (recall Section 2, Figure 7), whilst its orientation pattern agrees with
our convention.
\bigbreak

 $$ \picill4.4inby1.35in(D27) $$
 \begin{center}
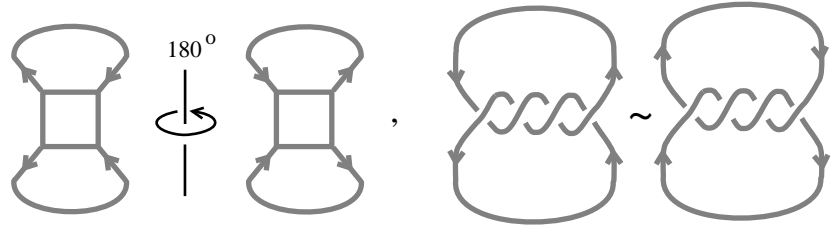
 {\bf Figure 27 -  Isotopic Oriented Rational Knots and Links }
\end{center}
 \vspace{3mm}

 Thus we reduce our analysis to two  basic
 types of orientation for the four end arcs of a rational tangle.
 We shall call an
 oriented rational tangle of {\it type I} if the {\it SW} arc is oriented
 upward and the {\it SE} arc
 is oriented downward, and of {\it type II} if the {\it SW} arc is oriented
 downward and the {\it SE} arc is oriented upward, see Figure~28. From the
 above remarks,  any tangle is of type
 I or type II. Two tangles are said to be {\it compatible} it they are both
 of type I or both of  type
 II and {\it incompatible} if they are of different types. In order to
classify oriented rational
 knots seen as numerator closures of oriented rational tangles, we will
always compare compatible
 rational rangles.  Note that if two oriented tangles are incompatible,
adding a single
 half twist at the bottom of one of them yields  a new pair of compatible
tangles, as Figure
 28 illustrates. Note also that adding such a twist, although it changes
the
tangle, it does not change the isotopy type of the numerator closure. Thus,
up to bottom
twists, we are always able to compare oriented rational tangles of the same
orientation type.
\bigbreak

$$ \picill3.7inby1.6in(D28) $$
\begin{center}
{\bf Figure 28 -  Compatible and Incompatible Orientations}
 \end{center}
 \vspace{3mm}

We shall now introduce  the notion of {\it connectivity} and we shall
relate it to orientation and the fraction of unoriented rational tangles.
We shall say that an unoriented rational tangle has {\it connectivity}
type $[0]$ if the {\it NW} end arc is connected to the {\it NE} end arc and the {\it
SW} end arc is connected to the  {\it SE} end arc.  Similarly, we say that
the tangle has {\it connectivity}  type $[+1]$ or type
$[\infty]$ if the end arc connections
are the same as in the tangles $[+1]$ and $[\infty]$ respectively.
The  basic connectivity patterns of rational tangles are exemplified by the
tangles $[0]$, $[\infty]$ and $[+1]$.  We can represent them iconically by the symbols
shown below.

 $$[0] = \mbox{\large $\asymp$}$$
 $$[\infty] = ><$$
 $$[+1] =\mbox{\large $\chi$}$$
\bigbreak

Note that connectivity type $[0]$ yields two-component rational
links, while  type
 $[+1]$ or $[\infty]$ yields one-component rational links. Also, adding a
bottom twist to a
rational tangle of  connectivity type $[0]$ will not change the connectivity
type of the
tangle, while  adding a bottom twist to a rational tangle of  connectivity
type $[\infty]$ will
switch the connectivity type to $[+1]$ and vice versa. While the
connectivity type of unoriented
rational tangles may be $[0]$, $[+1]$ or $[\infty],$ note that an oriented
rational
tangle of type I will have  connectivity  type $[0]$ or
$[\infty]$ and an oriented rational tangle of type II will have
connectivity  type $[0]$
or $[+1].$
\smallbreak

 Further, we need to keep an accounting of the  connectivity
 of rational  tangles in relation to the parity of the  numerators and
 denominators of their
 fractions.  We refer the reader to our paper \cite{KL2} for a full account.
\smallbreak

 We adopt the following notation:  $e$ stands for {\it even} and $o$ stands
 for {\it odd}. The {\it parity} of a fraction
 $p/q$ is defined to be the ratio of the parities ($e$ or $o$) of its
numerator and
 denominator $p$ and $q$. Thus the fraction $2/3$ is of parity $e/o.$
 The tangle $[0]$ has fraction $0 = 0/1,$ thus parity $e/o,$ the tangle
 $[\infty]$ has fraction $\infty = 1/0,$ thus parity $o/e,$ and
 the tangle $[+1]$ has fraction $1 = 1/1,$ thus parity $o/o.$ We then have
the following result.

 \begin{th}{ \ A rational tangle $T$ has  connectivity type
 $\mbox{\large $\asymp$}$ if and only if its fraction has parity $e/o$.
  $T$ has  connectivity type $><$ if  and only if its fraction has parity
 $o/e$. $T$ has  connectivity type $\mbox{\large $\chi$}$ if and only if
 its fraction has parity $o/o$. (Note that the formal fraction  of
$[\infty]$ itself is $1/0.$)
 Thus the link  $N(T)$ has two components if and only if $T$ has fraction
$F(T)$ of parity $e/o.$
  } \end{th}

 We will now proceed with sketching the proof of Theorem~3.
 We shall prove Schubert's oriented theorem by appealing to our previous
 work on the unoriented case and then analyzing how
 orientations and fractions are related. Our strategy is as follows:
Consider an oriented rational
 knot or link diagram $K$ in the form
 $N(T)$ where $T$ is a rational tangle in continued fraction form. Then any
other
 rational tangle  that
 closes to this knot $N(T)$ is available, up to bottom twists if necessary,
 as a cut from the given
 diagram. If two rational tangles close to give $K$ as an
 unoriented rational knot or link, then there are orientations on these
tangles,
 induced from $K$ so that the oriented tangles close
 to give $K$ as an oriented knot or link. The two tangles may or may not be
compatible.
 Thus, we must analyze when, comparing with the
 standard cut for the rational knot or link, another cut
 produces a compatible or incompatible rational tangle.
 However, assuming the top orientations are the same, we can replace one
 of the two incompatible tangles by the tangle obtained by adding a twist at
the
 bottom. {\it It is this possible twist difference that gives rise to the
change
 from modulus $p$ in the unoriented case to the modulus
 $2p$ in the oriented case.} We will now perform this
 analysis.  There
 are many interesting aspects to this analysis and we refer the reader to
 our paper
 \cite{KL2} for these details. Schubert \cite{Sch} proved his version
 of the oriented theorem
 by using the $2$-bridge representation of rational knots and  links, see
 also \cite{BZ}. We give a
 tangle-theoretic combinatorial proof based upon the combinatorics of the
 unoriented case.
\smallbreak

 The simplest instance of the classification of oriented
 rational knots is adding an
 {\it even number of twists} at the bottom of an oriented rational tangle
$T$, see
 Figure~28. We then obtain a
 compatible tangle $T*1/[2n],$ and $N(T*1/[2n]) \sim N(T).$ If now $F(T) =
 p/q$, then  $F(T*1/[2n]) = F(1/([2n] + 1/T)) = 1/(2n+1/F(T))=  p/(2np +
q).$
 Hence, if we set $2np + q = q'$
 we have $q\equiv q'(mod \, 2p),$ just as the oriented Schubert Theorem
 predicts. Note that reducing all
 possible bottom twists implies $|p|>|q|$ for both tangles, if the two
tangles that we compare each
time  are compatible or for only one, if they are incompatible.
 \smallbreak

We then have to  compare the special cut and the palindrome cut
with the
 standard cut.  In the oriented case the special cut is the easier to see
whilst
 the palindrome cut requires
 a more sophisticated analysis.  Figure 29 illustrates the general case of
 the special cut.  In order to
 understand Figure 29 it is necessary to also view Figure 22 for the details
of this cut.
\bigbreak

$$ \picill4.4inby1.7in(D29) $$
 \begin{center}
 {\bf Figure 29 - The Oriented Special Cut  }
\end{center}
 \vspace{3mm}

 Recall that if $S =
 [1]+T$ then the tangle of the special cut on the knot
 $N([1] + T)$ is the tangle $S' = [-1] -\frac{1}{T}.$ And if $F(T) = p/q$
 then $F([1]+T) =
 \frac{p+q}{q}$ and $F([-1] -\frac{1}{T}) = \frac{p+q}{-p}.$ Now, the point
 is that the orientations of
 the tangles $S$ and $S'$ are incompatible. Applying a $[+1]$ bottom twist
 to $S'$ yields
 $S'' = ([-1] -\frac{1}{T}) *[1]$, and we find that $F(S'') =
 \frac{p+q}{q}.$ Thus, the oriented
 rational tangles $S$ and $S''$ have the same fraction and by Theorem~1 and
 their compatibility they
 are oriented isotopic and the arithmetics of Theorem~3 is straightforward.
 \bigbreak

 We are left to examine the case of the palindrome cut. For this part of the proof, we refer the
reader to  our paper \cite{KL2}.

 \section{Strongly Invertible Links}

 An oriented knot or link is invertible if it is oriented
 isotopic to the link obtained from it by reversing the orientation of each component.
We have seen (see Figure 27) that  rational knots and
links are invertible.  A link $L$ of two components is said to
be {\it strongly invertible} if $L$ is ambient isotopic to itself with the orientation
of only one component reversed. In Figure 30 we illustrate the link
$L=N([[2],[1],[2]]).$ This is a strongly invertible link as is apparent by a $180^0$
vertical rotation.  This link is well-known as the Whitehead link, a link with linking
number zero. Note that since
$[[2], [1], [2]]$ has fraction equal to $1 + 1/(1 + 1/2) = 8/3$ this link is non-trivial
via the classification of rational knots and links. Note also that $3 \cdot 3 = 1 + 1
\cdot 8.$ 

\bigbreak

$$ \picill4.3inby2in(D30) $$

\begin{center}
{\bf Figure 30 - The Whitehead Link is Strongly Invertible } 
\end{center}
\vspace{3mm}

 
 In general we have the following. For our proof, see \cite{KL2}.

 \begin{th}\label{strongly} \ Let $L=N(T)$ be an oriented rational link with associated
 tangle fraction $F(T) = p/q$  of parity $e/o,$ with $p$ and $q$ relatively prime and
$|p|>|q|.$ Then $L$ is strongly invertible if and only if $q^{2} = 1 + up$ with $u$  an
odd integer. It follows that strongly invertible links are all numerators of rational
tangles  of the form $[[a_1], [a_2], \ldots, [a_k], [\alpha], [a_k], \ldots, [a_2],
[a_1]]$ for any integers $a_1, \ldots, a_k, \alpha.$
  \end{th}

\noindent See Figure 31 for another example of a strongly
invertible link.  In this case the link is $L = N([[3], [1], [1], [1], [3]])$ with $F(L)
= 40/11.$ Note that
$11^{2} =1 + 3\cdot 40,$ fitting the conclusion of Theorem \ref{strongly}.

\bigbreak

$$ \picill2.5inby1.6in(D31) $$

\begin{center}
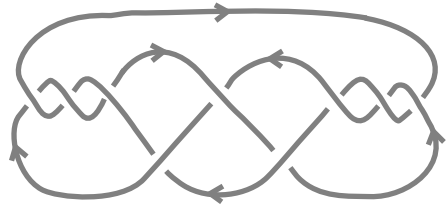
{\bf Figure 31 - An Example of a Strongly Invertible Link } 
\end{center}
\vspace{3mm}


 \section{Applications to the Topology of DNA}

 DNA supercoils, replicates and recombines with the help of certain
 enzymes.
 {\it Site-specific
 recombination} is one of the ways nature alters the genetic code of an
 organism, either by moving
 a block of DNA to another position on the molecule or by integrating a
 block
 of alien DNA into a
 host genome.   For a closed molecule of DNA a global picture of the
 recombination would be as shown in
 Figure 32, where double-stranded DNA is represented by a single line and
 the
 recombination sites are
 marked with points. This picture can be interpreted as $N(S + [0]) \longrightarrow
 N(S + [1]),$
 for $S = \frac{1}{[-3]}$ in this example. This operation can be repeated
 as in Figure 33. Note that
 the  $[0]-[\infty]$ interchange of Figure 10 can be seen as the first step
 of the process.
\bigbreak

$$ \picill4inby1.7in(D32) $$
\begin{center}
{\bf Figure 32 - Global Picture of Recombination}
 \end{center}
 \vspace{3mm}

 In this depiction of recombination, we have shown a local replacement of
 the tangle $[0]$
 by the tangle $[1]$ connoting a new cross-connection of the DNA strands.
 In general, it is not
 known without corroborating evidence just what the topological  geometry
 of the recombination
 replacement will be.  Even in the case of a single half-twist replacement
 such as $[1]$, it is
 certainly not obvious beforehand that the replacement will always be
 $[+1]$
 and not sometimes the
 reverse twist of $[-1].$ It was at the juncture raised by this question
 that
 a combination of
 topological methods in biology and a tangle model using knot theory
 developed by C.Ernst and D.W.~Sumners  resolved the issue in some specific cases. 
See \cite{ES}, \cite{Su} and references therein.
\bigbreak

$$ \picill4.3inby4.65in(D33) $$
\begin{center}
{\bf Figure 33 - Multiple Recombinations}
\end{center}
 \vspace{3mm}

 On the biological side, methods of protein coating developed by
 N. Cozzarelli,  S.J. Spengler and A. Stasiak et al. In \cite{CSS} it was
 made  possible for
 the first time to see knotted DNA in an electron micrograph with
 sufficient resolution to actually
 identify the topological type of these knots.  The protein coating
 technique made it possible to
 design an experiment involving successive DNA recombinations and to examine
 the topology of the products.  In
 \cite{CSS} the knotted DNA  produced by such successive recombinations was
 consistent with the
 hypothesis that all recombinations were of the type of a positive half
 twist as in $[+1].$ Then
 D.W.~Sumners  and C.~Ernst \cite{ES} proposed a {\em tangle model for
 successive DNA
 recombinations} and showed, in the case of the experiments in question,
 that there  was no other
 topological possibility for the recombination mechanism than the positive
 half twist $[+1].$ This
 constituted a unique use of topology as a theoretical
 underpinning for a problem in
 molecular biology.
 \bigbreak

 Here is a brief description of the tangle model for DNA
 recombination.  It is assumed
 that the initial state of the DNA is  described as the numerator closure
 $N(S)$ of a {\it substrate
 tangle} $S.$ The local geometry of the recombination is assumed to be
 described by the replacement
 of the tangle $[0]$ with a specific tangle $R.$ The results of the
 successive rounds of
 recombination are the knots and links
$$N(S + R) = K_1, \quad N(S + R + R) = K_2, \quad N(S + R + R + R) = K_3,
 \quad  \ldots$$

 \noindent  Knowing the knots $K_1, K_2, K_3, \ldots$ one would like to
 solve
 the above system of
 equations with the  tangles $S$ and $R$ as unknowns.  
\bigbreak

For such experiments
 Ernst and Sumners
 \cite{ES} used the classification of rational knots in the unoriented
 case,
 as well as results of
 Culler, Gordon, Luecke and Shalen \cite{CGLS} on Dehn surgery to prove that
 the
 solutions $S + nR$ must
 be {\it rational tangles}. These results of Culler, Gordon, Luecke and Shalen tell the topologist
under what circumstances a three-manifold with cyclic fundamental group must be a lens space. By showing 
when the $2$-fold branched covers of the DNA knots must be lens spaces, the recombination problems are reduced to the
consideration of  rational knots. This is a deep application of the three-manifold approach to rational knots and their
generalizations. 
\bigbreak

One can then apply the theorem on the
 classification of rational knots
 to deduce (in these instances) the  uniqueness of $S$ and $R$.  Note that,
 in these experiments, the
 substrate tangle $S$ was also pinpointed by the sequence of knots and
 links
 that resulted from the
 recombination.
 \bigbreak

 Here we shall solve tangle equations like the above under
 rationality  assumptions  on all tangles in question. This allows us to
 use
 only the mathematical
 techniques developed in this paper.
 We shall illustrate how a sequence of rational knots and links

 $$N(S+nR)=K_{n},  \ \ n=0,1,2,3, \ldots $$
\bigbreak

 \noindent with  $S$ and $R$ rational tangles, such that $R = [r], \ F(S) =
 \frac{p}{q}$
 and $p$, $q$, $r \in \ZZ$  ($p > 0$) {\em determines $\frac{p}{q}$ and
 $r$ uniquely}  if
 we know sufficiently many $K_{n}.$ We call this the ``DNA Knitting Machine
 Analysis".

 \begin{th} { Let a sequence $K_{n}$ of rational knots and links be defined
 by the equations  $K_{n} =
 N(S+nR)$ with specific integers $p$, $q$, $r$ ($p > 0$), where $R = [r], \
 F(S) = \frac{p}{q}.$  Then
 $\frac{p}{q}$ and $r$ are uniquely determined if one knows the topological
 type of the unoriented
 links $K_{0}, K_{1}, \ldots, K_{N}$ for any integer $N \geq |q| -
 \frac{p}{qr}.$ } \end{th}

 \noindent {\em Proof.} In this proof we shall write $N(\frac{p}{q} +nr)$ or
$N(\frac{p +qnr}{q})$
for $N(S+nR).$ We
 shall  also write  $K=K'$ to mean that $K$ and $K'$ are isotopic links.
 Moreover we shall say for a
 pair of reduced fractions $P/q$ and $P/q'$ that $q$ and $q'$ are {\it
 arithmetically related relative
 to $P$} if either $q\equiv q'(mod \, P)$   \, or \, $qq'\equiv 1(mod \,
 P).$  Suppose the integers $p, q, r$ give rise to the sequence of links
 $K_{0}, K_{1}, \ldots .$ Suppose there is some other triple of integers
$p', q', r'$ that give rise to the same sequence of links.
We will show uniqueness of $p, q, r$ under the conditions of the theorem.
We shall say ``the
 equality holds for $n$'' to mean that $N((p + qrn)/q) = N((p' +
 q'r'n)/q').$ We suppose that
 $K_{n} = N((p + qrn)/q)$ as in the hypothesis of the theorem, and suppose
 that there are $p'$, $q'$,
 $r'$ such that  for some $n$  (or a range of values of $n$ to be specified
 below) $K_{n} =  N((p' +
 q'r'n)/q').$
 \bigbreak

 If $n=0$ then we have $N(p/q) = N(p'/q').$ Hence by the
 classification theorem we know that
 $p=p'$  and that $q$ and $q'$ are arithmetically related. Note that the
 same
 argument shows that if
 the equality holds for any two  consecutive values of $n,$ then $p=p'.$
 Hence we shall assume
 henceforth that $p=p'.$ With this assumption in place, we see that if the
 equality holds for any $n
 \neq 0$  then $qr = q'r'.$ Hence we shall assume this as well from now on.
 \bigbreak

 If $|p + qrn|$ is sufficiently large, then the congruences for
 the arithmetical relation of
 $q$ and
 $q'$ must be {\em equalities over the integers}. Since $qq'= 1$ over the
 integers can hold only if
 $q=q'=1$ or $-1$ we see that it must be the case that $q = q'$ if the
 equality is to  hold for
 sufficiently large $n$. From this and the equation $qr = q'r'$ it follows
 that $r = r'.$
 It remains
 to determine a bound on $n$. In order to be sure that $|p + qrn|$ is
 sufficiently large, we need that
 $|qq'| \leq |p + qrn|.$  Since $q'r' = qr$, we also know that
 $|q'| \leq |qr|.$ Hence $n$ is sufficiently large if $|q^{2}r| \leq |p +
 qrn|.$
 \bigbreak

 If $qr > 0$ then,
 since $p > 0,$ we are asking that $|q^{2}r| \leq p + qrn.$  Hence
 $$n \geq (|q^{2}r| -p)/(qr) = |q| - (p/qr).$$ 
\bigbreak

If $qr < 0$ then for $n$
 large we will have $|p + qrn| = -p -qrn.$ Thus we want to solve
 $|q^{2}r| \leq -p - qrn$, whence
 $$n \geq (|q^{2}r| + p)/(-qr) = |q| - (p/qr).$$
\bigbreak

 Since these two cases exhaust the range of possibilities, this
 completes the proof of the
 theorem.
 $\hfill \Box$
 \bigbreak

 Here is  a special case of Theorem 7. See Figure 33.
 Suppose that we were given a sequence of knots and links  $K_{n}$ such
 that

 $$K_{n} = N( \frac{1}{[-3]} + [1]+ [1] +\ldots + [1] ) = N( \frac{1}{[-3]}
 +
 n\,[1]).$$

 \noindent We have  $F( \frac{1}{[-3]} + n\,[1]) = (3n-1)/3$ and we shall
 write $K_{n} =
 N([(3n-1)/3]).$ We are told that each of these rational knots is in fact
 the numerator closure of a rational tangle denoted

 $$[p/q] + n\,[r]$$

 \noindent for some rational number $p/q$ and some integer $r.$  That is,
 we are told that they
 come from a DNA knitting machine that is using rational tangle patterns.
 But we only know the
 knots  and the fact that they are indeed the closures for $p/q = -1/3$ and
 $r=1.$  By this analysis,
 the uniqueness is implied by the knots and links
 $\{ K_1, K_2, K_3, K_4 \}.$   This means that a DNA knitting machine
 $K_{n} = N(S + nR)$ that emits
 the four specific knots $K_n = N([(3n-1)/3])$ for $n=1,2,3,4$ must be of
 the form $S=1/[-3]$ and
 $R=[1]$. It was in this way (with a finite number of observations) that
 the structure of
 recombination in $T_{n3}$ resolvase was determined \cite{Su}.
 \bigbreak

 In this version of the tangle model for DNA recombination we have made a
 blanket assumption that
 the substrate tangle $S$ and the recombination tangle $R$ and all the
 tangles $S +nR$ were rational.
 Actually, if we assume that $S$ is rational and that $S+R$ is rational,
 then
 it follows
 that $R$ is an integer tangle.  Thus $S$ and $R$ neccesarily form a DNA
 knitting machine
 under these conditions. It is relatively natural to assume that $S$ is
 rational on the grounds of
 simplicity.  On the other hand it is not so obvious that the
 recombination
 tangle should be an
 integer. The fact that the products of the DNA recombination experiments
 yield rational knots and
 links, lends credence to the hypothesis of rational tangles and hence
 integral recombination
 tangles. But there certainly is a subtlety here, since we know that the
 numerator
 closure of the sum of two rational tangles is always a rational knot or
 link. In fact, it is here
 that some deeper topology shows that certain rational products from a
 generalized knitting machine
 of the form $K_n = N(S + nR)$  where $S$ and $R$ are arbitrary tangles
 will
 force the rationality of
 the tangles $S+nR$.  We refer the reader to \cite{ES}, \cite{EST},
\cite{ID} for the details of this approach.
 \bigbreak

\bigbreak

\noindent {\sc L.H.Kauffman: Department of Mathematics, Statistics and
 Computer Science, University
 of Illinois at Chicago, 851 South Morgan St., Chicago IL 60607-7045,
 U.S.A.}

 \vspace{.1in}
 \noindent {\sc S.Lambropoulou: Department of Mathematics, 
 National Technical University of Athens,
 Zografou campus, GR-157 80 Athens, Greece. }

\vspace{.1in}
\noindent {\sc E-mails:} \ {\tt kauffman@math.uic.edu  \ \ \ \ \ \ \ \ sofia@math.ntua.gr 

\noindent http://www.math.uic.edu/$\tilde{~}$kauffman/ \ \ \ \ http://users.ntua.gr/sofial}

 \end{document}